\documentclass{amsart}
\usepackage{latexsym,amsmath,amssymb,graphics,amscd}
\makeatletter
\@addtoreset{figure}{section}
\def\thefigure{\thesection.\@arabic\c@figure}
\def\fps@figure{h,t}
\@addtoreset{table}{bsection}
\def\thetable{\thesection.\@arabic\c@table}
\def\fps@table{h, t}
\@addtoreset{equation}{section}

\makeatother

\newcommand{\todo}[1]{\vspace{5 mm}\par \noindent
\marginpar{\textsc{ToDo}}
\framebox{\begin{minipage}[c]{0.95 \textwidth}
\tt #1 \end{minipage}}\vspace{5 mm}\par}
\newcommand{\bfi}{\bfseries\itshape}
\newtheorem{theorem}{Theorem}[section]

\newtheorem{corollary}[theorem]{Corollary}

\newtheorem{definition}[theorem]{Definition}

\newtheorem{lemma}[theorem]{Lemma}

\newtheorem{proposition}[theorem]{Proposition}
\newtheorem{remark}[theorem]{Remark}

\begin{document}
\title{Symmetry breaking for toral actions in simple mechanical
systems}
\author{Petre Birtea, Mircea Puta, Tudor S. Ratiu, R\u azvan
Tudoran}
\date{November 3, 2003}
\maketitle
\begin{abstract}
For simple mechanical systems, bifurcating branches of relative
equilibria with trivial symmetry from a given set of relative
equilibria with toral symmetry are found. Lyapunov stability
conditions along these branches are given.
\end{abstract}
\tableofcontents
\section{Introduction}
This paper investigates the problem of symmetry breaking in the
context of simple mechanical systems with compact symmetry Lie
group $G$. Let $\mathbb{T}$ be a maximal torus of $G$ whose Lie
algebra is denoted by $\mathfrak{t}$. Denote by $Q $ the
configuration space of the mechanical system. Assume that every
infinitesimal generator defined by an element of
$\mathfrak{t}$ evaluated at a symmetric configuration $q_e \in Q$
whose symmetry subgroup $G_{q_e}$ lies in $\mathbb{T}$ is a relative
equilibrium. The goal of this paper is to give sufficient conditions
capable to insure the existence of points in this set from which
branches of relative equilibria with trivial symmetry will emerge.
Sufficient Lyapunov stability conditions along  these branches
will be given if $G = \mathbb{T}$.
The strategy of the method can be roughly described
as follows. Denote by $\mathfrak{t} \cdot q_e$ the set of relative
equilibria described above. Take a regular element
$\mu\in \mathfrak{g}^{\ast }$ which happens to be the momentum
value of some relative equilibrium in $\mathfrak{t}\cdot q_{e}$. 
Choose a one parameter perturbation $\beta(\tau, \mu)  \in
\mathfrak{g}^\ast$ of $\mu$ that lies in
the set of regular points of $\mathfrak{g}^{\ast }$, for
small values of the parameter $\tau>0 $. Consider
the  $G_{q_{e}}$-representation on the tangent space
$T_{q_{e}}Q$.  Let $v_{q_e}$ be an element in the $\{e\}$-stratum 
of the representation  and also in the normal space to the tangent space at $q_{e}$ to the orbit
$G\cdot q_{e}$. Assume that its norm is small enough in order for
$v_{q_e}$ to lie in the open ball centered at the origin $0_{q_e}
\in T_{q_e}Q $ where the Riemannian exponential is a
diffeomorphism. The curve $\tau v_{q_e}$ projects by the
exponential map to a curve $q_e( \tau) $ in a neighborhood of
$q_e $ in $Q$ whose value at $\tau= 0$ is $q_e$. Note that the
isotropy subgroup at every point on this curve, except for $\tau=
0 $, is trivial. We shall search for relative equilibria in $TQ$
starting at points of
$\mathfrak{t}\cdot q_e $ such that their base curve in $Q $
equals $q_e(\tau)$ and their momentum values are $\beta(\tau,
\mu) $. To do this, we shall choose a curve $\xi(\tau,
v_{q_e}, \mu) \in \mathfrak{g} $ uniquely determined by $\beta(
\tau, \mu)$; as will be explained in the course of the
construction, $\xi(\tau, v_{q_e}, \mu)$ equals the value of
the inverse of the locked inertial tensor on $\beta(\tau, \mu)$
for $\tau \neq 0 $. If one can show that the limit of $\xi( \tau,
v_{q_e}, \mu) $ exists and belongs to $\mathfrak{t}$ for $\tau
\rightarrow  0$, then the infinitesimal generator of this value
evaluated at $q_e$ is automatically a relative equilibrium since 
it belongs to $\mathfrak{t}\cdot q_e$. It will be also shown that
the infinitesimal generators of $\xi(\tau, v_{q_e}, \mu)$
evaluated at $q_e(\tau)$ are relative equilibria. This produces a
branch of relative equilibria starting at this specific point in
$\mathfrak{t}\cdot q_e$ which has trivial isotropy for $\tau>0 $
and which depends smoothly on the additional parameter $\mu \in
\mathfrak{g}^\ast$.
In this method, there are two  key technical problems, namely, the
existence of the limit of $\xi( \tau,v_{q_e}, \mu)$ as $\tau
\rightarrow 0$ and the extension of the amended potential at
points with symmetry. The existence of the limit of 
$\xi( \tau,v_{q_e}, \mu)$ as $\tau \rightarrow 0$ will be
shown using the Lyapunov-Schmidt procedure. 
To extend the amended potential and its derivative at points with
symmetry, two auxiliary functions obtained by blow-up will
be introduced. The analysis breaks up in two problems on a space
orthogonal to the $G$-orbit.  
The present paper can be regarded as a sequel to the
work of Hern\'andez and Marsden \cite{hm}. The main difference
is that one single hypothesis from
\cite{hm} has been retained, namely that all points of
$\mathfrak{t} \cdot q_e $ are relative equilibria. We have also
eliminated a strong nondegeneracy assumption in \cite{hm}. But
the general principles of the strategy of the proof having to do
with a regularization of the amended potential at points with
symmetry, where it is not a priori defined, remains the same. In
a future paper we shall further modify this method to deal with
bifurcating branches of relative equilibria that have a given
isotropy, different from the trivial one, along the branch. 
The paper is organized as follows. In \S \ref{Lagrangian
mechanical systems} we quickly review the necessary material on
symmetric simple  mechanical systems and introduce the notations
and conventions for the entire paper. Relative equilibria and
their characterizations for general symmetric mechanical systems
and for simple ones in terms of the augmented and amended
potentials are recalled in \S \ref{Relative equilibria}. Section
\S \ref{Some basic results from the theory of group actions}
gives a brief summary of facts from the theory of proper group
actions needed in this paper. After these short introductory
sections, \S \ref{Regularization of the amended potential
criterion} presents the main bifurcation result of the
paper. The existence of branches of relative
equilibria starting at certain
points in $\mathfrak{t}\cdot q_e$, depending on several
parameters and having trivial symmetry off
$\mathfrak{t}\cdot q_e$, is proved in Theorem \ref{principala},
the main result of this paper. In \S
\ref{stability section}, using a result of Patrick \cite{patrick
thesis},  Lyapunov stability conditions for these branches are
given if the symmetry group is a torus.

\bigskip
\section{Lagrangian mechanical systems}
\label{Lagrangian mechanical systems}
This section summarizes the key facts from the theory of
Lagrangian systems with symmetry and sets the notations and
conventions to be used throughout this paper. The references for
this section are \cite{f of m}, \cite{lm}, \cite{marsden 92}, \cite{ims}.
\subsection{Lagrangian mechanical systems with symmetry }
Let $Q$ be a smooth manifold, the configuration space of a
mechanical system. The {\bfi fiber derivative\/} or {\bfi
Legendre transform\/} $\mathbb{F}L:TQ\rightarrow T^{\ast }Q$ of
$L$ is a vector bundle map covering the identity defined by
\[
\langle \mathbb{F}L(v_{q}),w_{q}\rangle
= \left.\frac{d}{dt}\right |_{t=0}L(v_{q}+tw_{q})
\]
for any $v_q, w_q \in TQ$. The {\bfi energy\/} of $L $ is
defined by $E(v_{q})= \langle \mathbb{F}L(v_{q}),v_{q} \rangle
-L(v_{q})$, $v_{q}\in T_{q}Q$. The pull back by $\mathbb{F}L $ of
the canonical one-- and two--forms of $T ^\ast Q $ give the
{\bfi Lagrangian one\/} and {\bfi two-forms\/} $\Theta _{L}$ and
$\Omega_L$ on $TQ$ respectively, that have thus the expressions
\[
\langle\Theta _{L}(v_{q}),\delta v_{q}\rangle=\langle\mathbb{F}L(v_{q}),T_{v_{q}}\pi
_{Q}(\delta v_{q})\rangle,\quad  v_{q}\in T_{q}Q,\quad  \delta v_{q}\in T_{v_{q}}TQ,
\qquad \Omega_{L}=-\mathbf{d}\Theta _{L},
\]
where $\pi_Q : TQ \rightarrow Q$ is the tangent bundle projection.
The Lagrangian $L $ is called {\bfi regular\/} if $\mathbb{F}L $
is a local diffeomorphism, which is equivalent  to $\Omega_L $
being a symplectic form on $TQ$. The Lagrangian $L $ is called
{\bfi hyperregular\/} if $\mathbb{F}L $ is a diffeomorphism and
hence a vector bundle isomorphism. The {\bfi Lagrangian vector
field\/} $X_E $ of $L$ is uniquely determined by the equality
\[
\Omega _{L}(v_{q})(X_{E}(v_{q}),w_{q})=\langle
\mathbf{d}E(v_{q}),w_{q}\rangle,
\quad \text{for} \quad v_{q}, \; w_{q}\in T_{q}Q.
\]
A {\bfi Lagrangian dynamical system\/}, or simply a {\bfi
Lagrangian system\/}, for $L$ is the dynamical system defined by
$X_{E}$, i.e., $\dot{v} =X_{E}(v)$. In standard coordinates $(q^i,
\dot{q}^i) $ the trajectories of $X_E $ are given by the second
order equations
\[
\frac{d}{dt}\frac{\partial L}{\partial \dot{q}^{i}}-\frac{\partial L}{
\partial q^{i}}=0,
\]
which are the classical the Euler-Lagrange equations.
\medskip
Let $\Psi: G \times Q \rightarrow Q$ be a smooth left Lie group
action on $Q$ and let
$L:TQ\rightarrow \mathbb{R}$ be a Lagrangian that is invariant
under the lifted action of $G$ to $TQ$. Denote by $\mathfrak{g}$
the Lie algebra of $G$. From the definition of the fiber
derivative it immediately follows that $\mathbb{F}L $ is
equivariant relative to the lifted $G $--actions  to $TQ$ and
$T^\ast Q$. The $G$-invariance of $L$ implies that
$X_{E}$ is $G$-equivariant, that is, $\Psi_g ^\ast X_E = X_E $
for any $g \in G $. The $G $--action on $TQ $ admits a momentum
map given by
\[
\langle\mathbf{J}_{L}\mathbf{(}v_{q}),\xi \rangle=\langle\mathbb{F}L(v_{q}),\xi
_{Q}(q)\rangle, \quad \text{for} \quad v_{q}\in T_{q}Q, \quad \xi
\in \mathfrak{g}.
\]
where $\xi _{Q}(q) : = d\exp (t\xi) \cdot q/dt|_{t=0}$ is the
{\bfi infinitesimal generator\/} of $\xi \in \mathfrak{g}$, where
$\mathfrak{g}$ denotes the Lie algebra of $G $. Recall that the
momentum map
$\mathbf{J}:T^{\ast}Q\rightarrow
\mathfrak{g}^{\ast }$ on $T^{\ast }Q$ is given by
\[
\langle\mathbf{J(}\alpha _{q}),\xi \rangle=\langle\alpha _{q},\xi _{Q}(q)\rangle,
\quad \text{for} \quad \alpha _{q}\in
T_{q}^{\ast }Q, \quad \xi \in \mathfrak{g}
\]
and hence $\mathbf{J}_{L}= \mathbf{J\circ }\mathbb{F}L$.
We shall denote by $g \cdot q: =\Psi(g,q)$ the action of the
element $g\in G $ on the point $q \in Q $. Similarly, the lifted
actions of $G $ on $TQ$ and $T^\ast Q $ are denoted by
\[
g \cdot v_q : = T_q \Psi_g(v_q) \quad \text{and} \quad
g \cdot \alpha_q : = T^\ast_{g\cdot q} \Psi_{g^{-1}}(\alpha_q)
\]
for $g\in G $, $v_q \in T_q Q $, and $\alpha_q \in T_q ^\ast Q $.
\subsection{Simple mechanical systems}
A {\bfi simple mechanical system} $(Q,\langle \!\langle \cdot ,\cdot \rangle \! \rangle
_{Q},V)$ consists of a Riemannian manifold $(Q,\langle \!\langle \cdot ,\cdot \rangle
\! \rangle _{Q})$ together with a potential function $V:Q\rightarrow \mathbb{R}$. These
elements define a Hamiltonian system on $(T^{\ast }Q,\omega )$ with Hamiltonian given
by $ H:T^{\ast }Q\rightarrow \mathbb{R}$, $H(\alpha _{q})= \frac{1}{2} \langle \!\langle
\alpha _{q},\alpha _{q}\rangle \! \rangle _{T^{\ast }Q}+V(q)$, where
$\alpha_q \in T_q ^\ast Q $ and $\langle \!\langle
\cdot ,\cdot \rangle
\! \rangle _{T^{\ast }Q}$ is the vector bundle metric on $T^{\ast
}Q$ induced by the Riemannian metric of $Q$. The Hamiltonian
vector field $X_H $ is uniquely given by the relation
$\mathbf{i} _{X_{H}}\omega =\mathbf{d}H$, where $\omega$ is the
canonical symplectic form on $ T^{\ast }Q$.
The dynamics of a simple mechanical system can also be described in terms of
Lagrangian mechanics, whose description takes place on $TQ$. The Lagrangian
for a simple mechanical system is given by $L:TQ\rightarrow \mathbb{R}$, $
L(v_{q})=\frac{1}{2}\langle\!\langle v_{q},v_{q}\rangle\!\rangle
_{Q}-V(q)$, where $v_q \in T_q Q $. The energy of $L $ is
$E(v_q) = \frac{1}{2}\langle\!\langle v_q, v_q \rangle\!\rangle +
V(q)$. Since the fiber derivative  for a simple mechanical system
is given by
$\langle\mathbb{F}L(v_{q}),w_{q}\rangle=\langle\!\langle
v_{q},w_{q}\rangle\!\rangle _{Q}$, or in local coordinates
$\mathbb{F}L\left(\dot{q}^{i}\frac{
\partial }{\partial q^{i}}\right)=g_{ij}\dot{q}^{j}dq^{i}$, where
$g_{ij}$ is the local expression for the metric on $Q$, it follows
that $L $ is hyperregular.
The relationship between the Hamiltonian and the Lagrangian
dynamics is the following: the vector bundle
isomorphism $\mathbb{F}L$ bijectively maps the trajectories of
$X_E$ to the trajectories of $X_H$, $(\mathbb{F}L)^\ast X_H = X_E
$, and the base integral curves of $X_E $ and $X_H $ coincide.
\subsection{Simple mechanical systems with symmetry}
\label{Simple mechanical systems with symmetry}
Let $G$ act on the configuration manifold $Q$ of a simple
mechanical system $(Q,\langle\!\langle \cdot ,\cdot
\rangle\!\rangle _{Q},V)$ by isometries.
The {\bfi locked inertia tensor\/} $\mathbb{I}:Q\rightarrow
\mathcal{L(} \mathfrak{g},\mathfrak{g}^{\ast })$, where
$\mathcal{L(} \mathfrak{g},\mathfrak{g}^{\ast })$ denotes the
vector space of linear maps from $\mathfrak{g}$ to
$\mathfrak{g}^\ast$,  is defined by
\[
\langle\mathbb{I(}q)\xi ,\eta \rangle=\langle\!\langle \xi
_{Q}(q),\eta _{Q}(q)\rangle\!\rangle _{Q}
\]
for any $q \in Q $ and any $\xi, \eta \in \mathfrak{g} $.
If the action is {\bfi locally free\/} at $q\in Q$, that is, the
isotropy subgroup $G_q $ is discrete, then
$\mathbb{I(}q)$  is an isomorphism and hence defines an inner
product on $\mathfrak{g}$. In general, the defining formula of
$\mathbb{I}(q)$ shows that $\ker \mathbb{I}(q) = \mathfrak{g}_q:
= \{\xi\in \mathfrak{g} \mid \xi_Q(q) = 0 \}$.
Suppose the action is locally free at every point $q \in Q $.
Then on can define the {\bfi mechanical connection\/}
$\mathcal{A} \in \Omega^1(Q; \mathfrak{g})$ by
\[
\mathcal{A}(q)(v_{q})=\mathbb{I(}q)^{-1}\mathbf{J}_{L}(v_{q}),
\quad v_{q}\in T_{q}Q.
\]
If the $G $--action is free and proper, so $Q \rightarrow Q/G $
is a $G $--principal bundle, then
$\mathcal{A}$ is a connection one--form on the principal bundle
$Q \rightarrow Q/G$, that is, it satisfies the following
properties:
\begin{itemize}
\item $\mathcal{A}(q) :T_q Q\rightarrow \mathfrak{g}$ is linear
and $G$-equivariant for every $q \in Q $, which means that
\[
\mathcal{A}(g \cdot q)(g\cdot v_q) =
\operatorname{Ad}_g[\mathcal{A}(q)(v_q)],
\]
for any  $v_q \in T_q Q $ and any $g \in G$, where
$\operatorname{Ad}$ denotes the adjoint representation of $G $ on
$\mathfrak{g}$;
\item $\mathcal{A}(q)(\xi _{Q}(q))=\xi $, for any $\xi \in
\mathfrak{g}$.
\end{itemize}
If $\mu\in \mathfrak{g}^\ast$ is given, we denote by
$\mathcal{A}_{\mu } \in \Omega^1(Q)$ the $\mu$--component of
$\mathcal{A}$, that is,  the one--form on
$Q$ defined by $\langle\mathcal{A}_{\mu }(q),v_{q}\rangle
=\langle\mu ,\mathcal{A}(q)(v_{q})\rangle$ for any $v_q \in T_qQ$.
The $G$-invariance of the metric and the relation
\[
(\operatorname{Ad}_{g}\xi
)_{Q}(q)=g\cdot \xi _{Q}(g^{-1}\cdot q),
\]
implies that
\begin{equation}
\label{equivariance of I}
\mathbb{I(}g\cdot q)=\operatorname{Ad}_{g^{-1}}^{\ast }\circ
\mathbb{I(}q)\circ \operatorname{Ad}_{g^{-1}}.
\end{equation}
We shall also need later the infinitesimal version of the above
identity
\begin{equation}
\label{infinitesimal equivariance of I}
T_q\mathbb{I} \left(\xi_{Q}(q)\right)
= -\operatorname{ad}_{\xi }^{\ast} \circ \mathbb{I}(q)
- \mathbb{I}(q) \circ \operatorname{ad}_{\xi},
\end{equation}
which implies
\begin{equation}
\label{useful identity for I}
\left\langle T_q \mathbb {I} ( \zeta_Q(q)) \xi, \eta
\right\rangle
= \mathbf{d}\langle\mathbb{I}(\cdot )\xi ,\eta
\rangle(q)\left(\zeta _{Q}(q)\right)
=\langle\mathbb {I}(q)[\xi,\zeta ],\eta
\rangle+\langle\mathbb{I(}q)\xi ,[\eta ,\zeta ]\rangle.
\end{equation}
for all $q\in Q$ and all $\xi, \eta, \zeta\in \mathfrak{g}$.

\section{Relative equilibria}
\label{Relative equilibria}
This section recalls the basic facts about relative equilibria
that will be needed in this paper. For proofs see \cite{f of m},
\cite{lm}, \cite{marsden 92}, \cite{ims}, \cite{slm}.
\subsection{Basic definitions and concepts}
Let $ \Psi: G \times Q \rightarrow Q$ be a left action of the Lie
group on the manifold
$Q$. A vector field $X:Q\rightarrow TQ$ is said to be $G$-{\bfi equivariant\/} if
\[
T_{q}\Psi _{g}(X(q)) = X(\Psi _{g}(q))\quad
\text{or,~equivalently,} \quad
\Psi_g^\ast X = X
\]
for all $q\in Q$ and $g\in G$. If $X$ is $G$-equivariant, then $G$ is said
to be a {\bfi symmetry group\/} of the dynamical system $\dot{q}=X(q)$.
A {\bfi relative equilibrium\/} of a $G$--equivariant vector field $X$ is a
point $q_{e}\in Q$ at which the value of $X$ coincides with the
infinitesimal generator of some element $\xi  \in \mathfrak{g}$,
usually called the {\bfi velocity\/} of $q_e $, i.e.,
\[
X(q_{e})=\xi _{Q}(q_{e}).
\]
A relative equilibrium $q_{e}$ is said to be {\bfi asymmetric\/}
if the isotropy subalgebra $\mathfrak{g}_{q_{e}} : = \{ \eta
\in \mathfrak{g} \mid \eta_Q(q_e) = 0 \} =\{0\}$, and {\bfi
symmetric} otherwise. Note that if
$q_{e}$ is a relative equilibrium with velocity $\xi \in
\mathfrak{g}$, then for any $g\in G$, $g\cdot q_{e}$ is a
relative equilibrium with velocity $\operatorname{Ad}_{g}\xi $.
The flow of an equivariant vector field induces a flow on the
quotient space. Thus, if the $G $--action is free and proper, a
relative equilibrium defines an equilibrium of the induced vector
field on the quotient space and conversely, any element in the
fiber over an equilibrium in the quotient space is a relative
equilibrium of the original system.
\subsection{Relative equilibria in Hamiltonian $G$-systems}
Given is a symplectic manifold $(P, \omega)$, a left Lie group
action of $G $ on $P $ that admits a momentum map $\mathbf{J}: P
\rightarrow \mathfrak{g}^\ast$, that is, $X_{\mathbf{J}^ \xi} =
\xi_P $, for any $\xi\in \mathfrak{g}$, where $\mathbf{J}^ \xi(p):
= \langle \mathbf{J}(p), \xi \rangle $, $p \in P$, is the
$\xi$--component of $\mathbf{J} $. We shall also assume
throughout this paper that the momentum map $\mathbf{J}$ is
equivariant, that is, $\mathbf{J}(g \cdot p ) =
\operatorname{Ad}^\ast_{g^{-1}} \mathbf{J} (p)$, for any $g \in G
$ and any $p \in P $. Given is also a $G $--invariant function
$H: P \rightarrow \mathbb{R}$. Noether's theorem states that the
$\mathbf{J} $ is conserved along the flow $F_t $ of the
Hamiltonian vector field $X_H $. In what follows we shall call
the quadruple $(Q,\omega ,H,\mathbf{J},G)$  a {\bfi Hamiltonian
$G$--system\/}. Consistent with the general definition presented
above, a point
$ p_{e}\in P$ is a {\bfi relative equilibrium\/} if
\[
X_{H}(p_{e})\in T_{p_{e}}(G\cdot p_{e}),
\]
where $G\cdot p_e : \{ g \cdot p_e \mid g \in G \} $ denotes the
$G $--orbit through $p_e $. Relative equilibria are
characterized in the following manner.
\begin{proposition}
\label{characterization of relative equilibria}
(\textbf{Characterization of relative equilibria}). Let $p_{e}\in
P$ and $ p_{e}(t)$ be the integral curve of $X_{H}$ with initial
condition $p_{e}(0)=p_{e}$. Let $\mu:= \mathbf{ J}(p_{e})$.
Then the following are equivalent:
\begin{enumerate}
\item[\textbf{(i)}] $p_{e}$ is a relative equilibrium.
\item[\textbf{(ii)}] There exists $\xi \in \mathfrak{g}$ such that
$p_{e}(t)=\exp (t\xi )\cdot p_{e}$.
\item[\textbf{(iii)}] There exists $\xi \in \mathfrak{g}$ such
that $p_{e}$ is a critical point of the {\bfi augmented
Hamiltonian\/}
\[
H_{\xi }(p):=H(p)-\langle\mathbf{J(}p)-\mu ,\xi \rangle.
\]
\end{enumerate}
\end{proposition}
Once we have a relative equilibrium, its entire $G$-orbit consists
of relative equilibria and the relation between the velocities of
the relative equilibria that are on the same $G$-orbit is given by
the adjoint action of $G$ on $\mathfrak{g}$.
\begin{proposition}
\label{orbita}
With the notations  of the previous proposition, let $p_{e}$ be a
relative equilibrium with velocity $\xi$. Then
\begin{enumerate}
\item[\textbf{(i)}] for any
$g\in G$, $g\cdot q_{e}$ is also a relative
equilibrium whose velocity is $\operatorname{Ad}_g \xi$;
\item[\textbf{(ii)}] $\xi(q_{e})\in \mathfrak{g}_\mu : = \{\eta
\in \mathfrak{g}\mid \operatorname{ad}^\ast_\eta \mu = 0 \}$, the
coadjoint isotropy subalgebra at $\mu\in \mathfrak{g}^\ast$, i.e.,
$\operatorname{Ad}_{\exp t\xi}^{\ast}\mu=\mu$ for any $t \in
\mathbb{R}$.
\end{enumerate}
\end{proposition}
\subsection{Relative equilibria in simple mechanical $G$-systems}
In the case of simple mechanical $G$-systems, the
characterization $\mathbf{(iii)}$ in Proposition
\ref{characterization of relative equilibria} can be simplified in
such way that the search of relative equilibria reduces to the
search of critical points of a real valued function on
$Q$.
Depending on whether one \ keeps track of the velocity or
the momentum of a relative equilibrium, this simplification
yields the \textit{augmented} or the \textit{amended} potential
criterion, which we introduce in what follows.
Let $(Q,\langle\!\langle \cdot ,\cdot \rangle\!\rangle _{Q},V,G)$
be a simple mechanical $G$--system.
\begin{itemize}
\item For $\xi \in \mathfrak{g}$, the {\bfi augmented potential\/}
$V_{\xi }:Q\rightarrow \mathbb{R}$ is defined by $V_{\xi
}(q):=V(q)-\frac{1}{2}\langle
\mathbb{I}(q)\xi ,\xi \rangle$.
\item For $\mu \in \mathfrak{g}^{\ast }$, the {\bfi amened
potential\/} $V_{\mu }:Q\rightarrow \mathbb{R}$ is defined by
$V_{\mu }(q):=V(q)+\frac{1}{2}\langle
\mathbb{\mu },\mathbb{I}(q)^{-1}\mu \rangle$.
\end{itemize}
Note that the amended potential is defined at $q\in Q$ only if $q$
in an asymmetric point.
There is an alternate expression for the amended potential,
namely,  $V_{\mu }(q)=(H\circ \mathcal{A}_{\mu })(q)$.
\begin{proposition}
(\textbf{Augmented potential criterion}).
\label{augmented potential criterion}
A point
$(q_{e},p_{e})\in T^{\ast}Q$ is a relative equilibrium if and only
if there exists a $\xi \in
\mathfrak{g}$ such that:
\begin{enumerate}
\item[\textbf{(i)}] $p_{e}=\mathbb{F}L(\xi _{Q}(q_{e}))$ and
\item[\textbf{(ii)}] $q_{e}$ is a critical point of $V_{\xi }$.
\end{enumerate}
\end{proposition}
\begin{proposition}
(\textbf{Amended potential criterion}).
\label{amended potential criterion}
A point $(q_{e},p_{e})\in
T^{\ast }Q$ is a relative equilibrium if and only if there exists
a $\mu \in \mathfrak{g} ^{\ast }$ such that:
\begin{enumerate}
\item[\textbf{(i)}] $p_{e}=\mathcal{A}_{\mu }(q_{e})$ and
\item[\textbf{(ii)}] $q_{e}$ is a critical point of $V_{\mu }$.
\end{enumerate}
\end{proposition}

\section{Some basic results from the theory of Lie group actions}
\label{Some basic results from the theory of group actions}
We shall need a few fundamental results form the theory of group
actions which we now review. For proofs and further information
see \cite{br}, \cite{dk}, \cite{kawakubo}, \cite{or}.
\subsection{Maximal tori}
Let $V$ be a representation space of a compact Lie group $G$.
A point $v\in V$ is {\bfi regular\/} if there is no $G$--orbit in
$V$ whose dimension is strictly greater than the dimension of the
$G$--orbit through $v$. The set of regular points, denoted
$V_{reg}$,  is open and dense in $V$.
In particular, $\mathfrak{g}_{reg}$ and
$\mathfrak{g}_{reg}^{\ast}$, denote the set of regular points in
$\mathfrak{g}$ and $\mathfrak{g}^{\ast }$ with respect to adjoint
and coadjoint representation, respectively.
A subgroup of a Lie group is said to be a {\bfi torus\/} if it is
isomorphic to $S^{1}\times \cdot \cdot \cdot \times S^{1}$.
Every Abelian subgroup of a compact connected Lie group is
isomorphic to a torus.
A subgroup of a Lie group is said to be a {\bfi maximal torus\/}
if it is a torus that is not properly contained in some other
torus.
Every $\xi \in \mathfrak{g}$ belongs to at least one maximal
Abelian subalgebra and every $\xi \in \mathfrak{g\cap g}_{reg}$
belongs to exactly one such maximal Abelian subalgebra.
Every maximal Abelian subalgebra is the Lie algebra of some
maximal torus  in $G$.
Let $\mathfrak{t}$ be the maximal Abelian subalgebra corresponding
to a maximal torus $T$. Then for any $\xi \in \mathfrak{t\cap
g}_{reg}$, we have that $G_{\xi }=T$. The space
$[\mathfrak{g},\mathfrak{t}]$ is the orthogonal complement to
$\mathfrak{t}$ in $\mathfrak{g}$ with respect to any $G$--invariant
inner product on $\mathfrak{g}$. Such an inner product exists by
compactness of $G$ by simply averaging any inner product on
$\mathfrak{g}$. Therefore, we have
$\mathfrak{ g} =\mathfrak{t}\oplus [\mathfrak{g},\mathfrak{t}]$.
Let $[\mathfrak{g},\mathfrak{t}]^{\circ }$ the annihilator of
$[\mathfrak{g},\mathfrak{t}]$. Then $G_{\mu }=T$ for every $\mu
\in [\mathfrak{g},\mathfrak{t} ]^{\circ }\cap
\mathfrak{g}_{reg}^{\ast }$. Since
$[\mathfrak{g},\mathfrak{t}]^{\circ }\cap
\mathfrak{g}_{reg}^{\ast }$ is dense in
$[\mathfrak{g},\mathfrak{t}]^{\circ }$, it follows that
$T\subset G_{\mu }$ for every $\mu \in
[\mathfrak{g},\mathfrak{t}]^{\circ }$.
\subsection{Twisted products}
Let $G$ be a Lie group and $H\subset G$ be a Lie subgroup. Suppose
that $H$ acts on the left on a manifold $A$. The {\bfi twisted
action\/} of $H$ on the product $G\times A$ is defined by
\[
h\cdot (g,a)=(gh,h^{-1}\cdot a), \quad h\in H, \quad g\in G,
\quad a\in A.
\]
Note that this action is free and proper by the freeness and properness of
the action on the $G$--factor. The {\bfi twisted product\/}
$G\times _{H}A$ is defined as the orbit space $(G\times A)/H$
of the twisted action. The elements of $G\times
_{H}A$ will be denoted by $[g,a],$ $g\in G,$
$a\in A$.
The twisted product $G\times _{H}A$ is a $G$--space relative to
the left action defined by $g^{\prime }\cdot \lbrack
g,a]=[g^{\prime }g,a]$. Also, the action of $H$ on $A$ is proper
if and only if the $G$--action on $G\times _{H}A$ is proper.
The isotropy subgroups of the $G$--action on the twisted product
$G\times _{H}A$ satisfy
\[
G_{[g,a]}=gH_{a}g^{-1}, \quad g\in G, \quad a\in A.
\]
\subsection{Slices}
Throughout this paragraph it will be assumed that $\Psi : G
\times Q \rightarrow Q $ is a left proper action of the Lie group
$G$ on the manifold $Q$. This action
will not be assumed to be free, in general. For $q\in Q$ we will
denote by $H:=G_{q} := \{ g \in G \mid g \cdot q = q \}$ the
isotropy subgroup of the action $\Psi$ at $q $.
We shall introduce also the following convenient notation: if $K
\subset G $ is a Lie subgroup of $G $ (possibly equal to $G $),
$\mathfrak{k}$ is its Lie algebra,  and $q\in Q$, then
$\mathfrak{k}\cdot q := \{ \eta_Q(q) \mid \eta \in
\mathfrak{k} \}$ is the tangent space to the orbit $K\cdot q $ at
$q $.
A {\bfi tube\/} around the orbit
$G\cdot q$ is a $G$-equivariant diffeomorphism $\varphi :G\times
_{H}A\rightarrow U$, where $U$ is a $G$-invariant neighborhood of
$G\cdot q$ and $A$ is some manifold on which $H$ acts. Note that
the $G$-action on the twisted product $ G\times _{H}A$ is proper
since the isotropy subgroup $H$ is compact and, consequently, its
action on $A$ is proper. Hence the $G$-action on $G\times _{H}A$
is proper.
Let $S$ be a submanifold of $Q$ such that $q\in S$ and $H\cdot
S=S$. We say that $S$ is a {\bfi slice\/} at $q$ if the map
\[
\varphi :G\times _{H}S\rightarrow U
\]
\[
\lbrack g,s]\mapsto g\cdot s
\]
is a tube about $G\cdot q$, for some $G$--invariant open
neighborhood of $G\cdot q$. Notice that if $S$ is a slice at
$q$ then $g\cdot S$ is a slice at the point $g\cdot q$.
The following statements are equivalent:
\begin{enumerate}
\item[\textbf{(i)}] There is a tube $\varphi :G\times
_{H}A\rightarrow U$ about $G\cdot q$ such that $\varphi
([e,A])=S$.
\item[\textbf{(ii)}] $S$ is a slice at $q$.
\item[\textbf{(iii)}] The submanifold $S$ satisfies the following
properties:
\begin{enumerate}
\item[\textbf{(a)}] The set $G\cdot S$ is an open neighborhood of
the orbit $G\cdot q$ and
$S$ is closed in $G\cdot S$.
\item[\textbf{(b)}] For any $s\in S$ we have
$T_{s}Q=\mathfrak{g}\cdot s+T_{s}S$. Moreover, $\mathfrak{g}\cdot
s\cap T_{s}S=\mathfrak{h}\cdot s$, where $\mathfrak{h}$ is the
Lie algebra of $H $. In particular
$T_{q}Q=\mathfrak{g}\cdot q\oplus T_{q}S$.
\item[\textbf{(c)}] $S$ is $H$-invariant. Moreover, if $s\in S$
and $g\in G$ are such that
$g\cdot s\in S$, then $g\in H$.
\item[\textbf{(d)}] Let $\sigma :U\subset G/H\rightarrow G$ be a
local section of the submersion $G\rightarrow G/H$. Then the map
$F:U\times S\rightarrow Q$ given by $F(u,s):=\sigma (u)\cdot s$
is a diffeomorphism onto an open set of $Q$.
\end{enumerate}
\item[\textbf{(iv)}] $G\cdot S$ is an open neighborhood of $G\cdot
q$ and there is an equivariant smooth retraction
\[
r:G\cdot S\rightarrow G\cdot q
\]
of the injection $G\cdot q\hookrightarrow G\cdot S$ such that $r^{-1}(q)=S$.
\end{enumerate}
\begin{theorem}
\textbf{(Slice Theorem)} Let $Q$ be a manifold and $G$ be a Lie
group acting properly on $Q$ at the point $q\in Q$. Then, there
exists a slice for the $G$--action at $q$.
\end{theorem}
\begin{theorem}
\label{tube theorem}
\textbf{(Tube Theorem)} Let $Q$ be a manifold and $G$ be a Lie group acting
properly on $Q$ at the point $q\in Q$, $H:=G_{q}$. There exists a tube $
\varphi :G\times _{H}B\rightarrow U$ about $G\cdot q$ such that $\varphi
([e,0])=q$, $\varphi ([e,B])=:S$ is a slice at $q$; $B$ is an
open  $H$--invariant neighborhood of $0$ in the vector space
$T_{q}Q/T_{q}(G\cdot q) $, on which $H$ acts linearly by $h\cdot
(v_{q}+T_{q}(G\cdot q)):=T_{q}\Psi _{h}(v_{q})+T_{q}(G\cdot
q)$.
\end{theorem}
If $Q$ is a Riemannian manifold then $B$ can be chosen to be a
$G_{q}$--invariant neighborhood of $0$ in $(\mathfrak{g}\cdot
q)^{\perp }$, the orthogonal complement to $\mathfrak{g}\cdot q$
in $T_{q}Q$. In this case $ U=G\cdot \operatorname{Exp}_{q}(B)$,
where $\operatorname{Exp}_{q}: T_q Q \rightarrow Q $ is the
Riemannian exponential map.
\subsection{Type submanifolds and fixed point subspaces}
Let $G$ be a Lie group acting on a manifold $Q$. Let $H$ be a closed
subgroup of $G$. We define the following subsets of $Q$ :
\begin{eqnarray*}
Q_{(H)} &=&\{q\in Q\mid G_{q}=gHg^{-1},g\in G\}, \\
Q^{H} &=&\{q\in Q\mid H\subset G_{q}\}, \\
Q_{H} &=&\{q\in Q\mid H=G_{q}\}.
\end{eqnarray*}
All these sets are submanifolds of $Q $.
The set $Q_{(H)}$ is called the $(H)$--{\bfi orbit type
submanifold\/},
$ Q_{H}$ is the $H$--{\bfi isotropy type submanifold\/}, and
$Q^{H}$ is the $H$--{\bfi fixed point submanifold\/}. We
will collectively call these subsets the {\bfi type
submanifolds}. We have:
\begin{itemize}
\item  $Q^{H}$ is closed in $Q$;
\item $Q_{(H)}=G\cdot Q_{H}$;
\item $Q_H $ is open in $Q^H $.
\item the tangent space at $q \in Q_H $ to $Q_H $ equals
\[
T_{q}Q_{H}=\{v_{q}\in T_{q}Q\mid T_{q}\Psi_{h}(v_{q})
=v_{q},\, \forall h\in H\}=(T_q Q)^H=T_q Q^H;
\]
\item $T_{q}(G\cdot q)\cap (T_{q}Q)^{H}=T_{q}(N(H)\cdot q)$,
where $N(H) $ is the normalizer of $H $ in $G $;
\item if $H$ is compact then $Q_{H}=Q^{H}\cap Q_{(H)}$ and $Q_H $
is closed in $Q_{(H)}$.
\end{itemize}
If $Q$ is a vector space on which $H$ acts linearly, the set
$Q^{H}$ is found in the physics literature under the names of
{\bfi space of singlets\/} or {\bfi space of invariant vectors\/}.
\begin{theorem}
\textbf{(The stratification theorem)}. Let $Q$ be a smooth manifold and $G$
be a Lie group acting properly on it. The connected components of the orbit
type manifolds $Q_{(H)}$ and their projections onto orbit space $Q_{(H)}/G$
constitute a Whitney stratification of $Q$ and $Q/G$, respectively. This
stratification of $Q/G$ is minimal among all Whitney stratifications of $Q/G$.
\end{theorem}
The proof of this result, that can be found in \cite{dk} or
\cite{pflaum}, is based on the Slice Theorem and on a series
of extremely important properties of the orbit type manifolds
decomposition that we enumerate in what follows. We start by
recalling that the set of conjugacy classes of subgroups of a Lie
group $G$ admits a partial order by defining $ (K)\preceq (H)$ if
and only if $H$ is conjugate to a subgroup of $K$. Also, a point
$q\in Q$ in a proper $G$--space $Q$ (or its corresponding
$G$--orbit, $ G\cdot q$) is called {\bfi principal\/} if its
corresponding local orbit type manifold is open in $Q$. The orbit
$G\cdot q$ is called
{\bfi regular\/} if the dimension of the orbits nearby coincides
with the dimension of $ G\cdot q$. The set of principal and
regular orbits will be denoted by $ Q_{princ}/G$ and $Q_{reg}/G$,
respectively. Using this notation we have:
\begin{itemize}
\item  For any $q\in Q$ there exists an neighborhood $U$ of $q$
that intersects only finitely many connected components of
finitely many orbit type manifolds. If $Q$ is compact or a linear
space where $G$ acts linearly, then the $G$--action on $Q$ has
only finitely many distinct connected components of orbit type
manifolds.
\item  For any $q\in Q$ there exists an open neighborhood $U$ of
$q$ such that $(G_{q})\preceq (G_{x})$, for all $x\in U$. In
particular, this implies that $\dim G\cdot q\leq \dim G\cdot x$,
for all
$x\in U$.
\item  {\bfi Principal Orbit Theorem}: For every connected
component
$ Q^{0}$ of $Q$ the subset $Q_{reg}\cap Q^{0}$ is connected,
open, and dense in $Q^{0}$. Each connected component $(Q/G)^{0}$
of $Q/G$ contains only one principal orbit type, which is
connected open and dense in $(Q/G)^0 $.
\end{itemize}

\section{Regularization of the amended potential criterion}
\label{Regularization of the amended potential criterion}
In this section we shall follow the strategy in \cite{hm} to give
sufficient criteria for finding relative equilibria emanating
from a given one and to find a method that distinguishes
between the distinct branches. The criterion will involve a
certain regularization of the amended potential. The main
difference with \cite{hm} is that all hypotheses but one have been
eliminated and we work with a general torus and not just a
circle. The conventions, notations, and method of proof are those
in \cite{hm}.
\subsection{The bifurcation problem}
Let $(Q,\langle\!\langle \cdot ,\cdot \rangle\!\rangle _{Q},V,G)$
be a simple mechanical
$G$-system, with $G$ a compact Lie group with the Lie algebra
$\mathfrak{g}$. Recall that the left  $G $--action $\Psi:G \times
Q \rightarrow Q $ is by isometries and that the potential $V:Q
\rightarrow \mathbb{R}$ is $G$--invariant. Let $ q_{e}\in Q$ be a
symmetric point whose isotropy group $G_{q_{e}}\subset
\mathbb{T}$ is contained in a maximal torus $\mathbb{T}$ of $G$.
Denote by $\mathfrak{t}
\subset \mathfrak{g}$ the Lie algebra of $\mathbb{T}$; thus
$\mathfrak{t}$ is a maximal Abelian Lie subalgebra of
$\mathfrak{g}$. Throughout this section we shall make the
following hypothesis:
\smallskip
\begin{center}
\textbf{(H)} \textit{ every} $v_{q_{e}}\in \mathfrak{t}\cdot
q_{e}$ \textit{is~a~relative~equilibrium}.
\end{center}
\smallskip
The following result was communicated to us by J. Montaldi.
\begin{proposition}
\label{montaldi}
In the context above we have that:
\[
\begin{array}{cc}
(i) & \mathbf{d}V(q_{e})=0 \\
(ii) & \mathbb{I(}q_{e}\mathbb{)}\mathfrak{t}\subseteq
 [\mathfrak{g},
\mathfrak{t}]^{\circ }.
\end{array}
\]
\end{proposition}
\begin{proof}
$(i)$ Because all the elements in $\mathfrak{t}\cdot q_{e}$ are relative
equilibria, we have by the augmented potential criterion
$\mathbf{d}V_{\xi }(q_{e})=0$, for any $\xi \in \mathfrak{t}$.
Consequently for $\xi =0$ we will obtain
$0=\mathbf{d}V_{0}(q_{e})=\mathbf{d}V(q_{e})$.
$(ii)$ Substituting in the relation \eqref{useful identity for I}, $q$ by
$q_{e}$ and setting  $\eta =\xi
\in \mathfrak{t}$ we obtain:
\[
\mathbf{d}\langle\mathbb{I(}\cdot \mathbb{)}\xi ,\xi \rangle(q_{e})
(\zeta _{Q}(q_{e})) =\langle
\mathbb{I(}q_{e}\mathbb{)[}\xi ,\zeta ],\xi \rangle
+\langle\mathbb{I(}q_{e}\mathbb{)}\xi
,[\xi ,\zeta ]\rangle = 2\langle \mathbb{I(}q_{e}\mathbb{)}\xi
,[\xi ,\zeta ]\rangle
\]
for any $\xi \in \mathfrak{t}$ and $\zeta \in \mathfrak{g}$.
The  augmented potential criterion yields
\[
0 = \mathbf{d}
V_{\xi}(q_{e})
=\mathbf{d}V(q_{e})-\frac{1}{2}\mathbf{d}\langle\mathbb{I(}\cdot
\mathbb{)}\xi ,\xi
\rangle(q_{e}).
\]
Since $\mathbf{d}V(q_{e})=0$ by (i), this implies
$\mathbf{d}\langle\mathbb{I(}\cdot
\mathbb{)}\xi ,\xi \rangle(q_{e})=0$ and consequently
$\langle\mathbb{I(}q_{e}\mathbb{)}\xi , [\xi ,\zeta ]\rangle=0$,
for any $\xi \in \mathfrak{t}$ and $\zeta \in \mathfrak{g}$. So we have the inclusion
\[
\mathbb{I(}q_{e}\mathbb{)}\xi \subseteq [\mathfrak{g},\xi
]^{\circ }.
\]
Now we will prove that $[\mathfrak{g},\xi ]^{\circ }=[\mathfrak{g},\mathfrak{
t}]^{\circ }$ for regular elements $\xi \in \mathfrak{t}$. For this it is enough
to prove that $[\xi ,\mathfrak{g}]=[\mathfrak{t},
\mathfrak{g}]$ for regular elements $\xi \in \mathfrak{t}$.
It is obvious that $[\xi ,\mathfrak{g}]\subseteq \lbrack \mathfrak{t},
\mathfrak{g}]$. Equality will follow by showing  that  both spaces
have the same dimension. To do this, let
$F_{\xi }:\mathfrak{g}\rightarrow
 \mathfrak{g}$, $F_{\xi }(\eta):=\operatorname{ad}_{\xi }\eta $,
which is obviously a linear  map whose image and kernel are
$\operatorname{Im} (F_{\xi })=[\xi ,\mathfrak{g}]$ and
$\ker (F_{\xi })=\mathfrak{g} _{\xi }$. Because $\xi \in
\mathfrak{t}$ is a regular element we have that $\mathfrak{g}_{\xi
}=\mathfrak{t}$ and so $\ker (F_{\xi })=\mathfrak{t}$. Thus
$\dim (\mathfrak{g})=\dim (
\mathfrak{t})+\dim ([\xi ,\mathfrak{g}])$ and so using the fact that $\dim (
\mathfrak{g})=\dim (\mathfrak{t})+\dim ([\mathfrak{t},\mathfrak{g}])$ (since
$\mathfrak{g}=\mathfrak{t}\oplus [\mathfrak{g},\mathfrak{t}]$,
$\mathfrak{g}$ being a compact Lie algebra), we obtain the equality
$\dim ([\xi ,\mathfrak{g}])=\dim ([\mathfrak{t},\mathfrak{g}])$.
Therefore,  $[\xi,\mathfrak{g}]=[\mathfrak{t},\mathfrak{g}]$ for
any regular element $\xi \in \mathfrak{t}$.
Summarizing, we proved
\[
\mathbb{I(}q_{e}\mathbb{)}\xi \subseteq
[\mathfrak{g},\mathfrak{t} ]^{\circ },
\]
for any regular element $\xi \in \mathfrak{t}$.
The continuity of $\mathbb{I(}q_{e}\mathbb{)}$, the closedness of
$[\mathfrak{g},\mathfrak{t}]^{\circ }$, and that fact that the
regular elements $\xi \in \mathfrak{t}$ form a dense subset of
$\mathfrak{t}$, implies that
\[
\mathbb{I(}q_{e}\mathbb{)}\xi \subseteq
[\mathfrak{g},\mathfrak{t} ]^{\circ },
\]
for any $\xi \in \mathfrak{t}$ and hence $\mathbb{I(}q_{e}\mathbb{)}
\mathfrak{t}\subseteq [\mathfrak{g},\mathfrak{t}]^{\circ }$.
\end{proof}
\begin{lemma}
\label{same isotropy}
For each $v_{q_{e}}\in \mathfrak{t}\cdot q_{e}$ we have
$G_{v_{q_{e}}}=G_{_{q_{e}}}$.
\end{lemma}
\begin{proof}
The inclusion $G_{v_{q_{e}}}\subseteq G_{_{q_{e}}}$ is obviously
true, so it will be enough to prove that $G_{v_{q_{e}}}\supseteq
G_{_{q_{e}}}$. To see this, let $g\in G_{_{q_{e}}}$ and
$v_{q_{e}}=\xi _{Q}(q_{e})\in
\mathfrak{t}\cdot q_{e}$, with $\xi\in \mathfrak{t}$. Then, since
$G_{q_e} $ is Abelian, we get
\begin{align*}
T_{q_{e}}\Psi _{g}\left(v_{q_{e}}\right)
&=T_{q_{e}}\Psi _{g}\left(\xi_{Q}(q_{e})\right)
=T_{q_{e}}\Psi _{g}\left(\left.\frac{d}{dt}\right|_{t=0}\Psi
_{\exp (t \xi )}(q_{e})\right) \\
&=\left.\frac{d}{dt}\right|_{t=0}\left(\Psi _{g}\circ \Psi _{\exp
(t\xi )}\right)(q_{e})=
\left.\frac{d}{dt}\right| _{t=0}(\Psi _{\exp (t \xi )}\circ \Psi
_{g})(q_{e}) \\ &=\left.\frac{d}{dt}\right|_{t=0}\Psi _{\exp
(t\xi )}(q_{e})= \xi_Q(q_e) = v_{q_{e}},
\end{align*}
that is, $g\cdot v_{q_{e}} = v_{q_{e}} $, as required.
\end{proof}
The bifurcation problem for relative equilibria on $TQ$ can be
regarded as a bifurcation problem on the space $Q\times\mathfrak{g}
^{\ast}$ as the following shows.
\begin{proposition}
\label{map f}
The map $f:TQ\rightarrow Q\times \mathfrak{g}^{\ast }$ given by $
v_{q}\mapsto (q,\mathbf{J}_{L}(v_{q}))$ restricted to the set
of relative equilibria is one to one and onto its image.
\end{proposition}
\begin{proof}
The only thing to be proved is that the map is injective. To see
this, let  $(q_{1},(\xi_{1})_{Q}(q_{1}))$ and
$(q_{2},(\xi_{2})_{Q}(q_{2}))$  be two relative equilibria such
that $f(q_{1},(\xi_{1})_{Q}(q_{1}))=
f(q_{2},(\xi_{2})_{Q}(q_{2}))$.  Then $q_{1}=q_{2}=:q$ and
$\mathbf{J}_{L}(q,(\xi _{1}-\xi _{2})_{Q}(q))=
\mathbb{I}(q)(\xi_{1}-\xi_{2})=0$ which shows that $\xi_{1}-\xi_{2}\in \ker
\mathbb{I}(q)=\mathfrak{g}_{q}$ and hence $(\xi_{1})_
{Q}(q)=(\xi_{2})_{Q}(q)$.
\end{proof}
We can thus change the problem: instead of searching for relative
equilibria of the simple mechanical system in $TQ$, we shall set
up a bifurcation problem on $Q \times \mathfrak{g}^\ast$ such
that the image of the relative equilibria by the map $f $ is
precisely the bifurcating set.
To do this, we begin with some geometric considerations.  We
construct a $G$-invariant tubular neighborhood of the orbit
$G\cdot q_{e}$ such that the isotropy group of every point in
this neighborhood is a subgroup of $G_{q_{e}}$. This follows from
the Tube Theorem \ref{tube theorem}. Indeed, let
$B\subset (\mathfrak{g}\cdot q_{e})^{\perp }$ be a
$G_{q_{e}}$-invariant open neighborhood of
$0_{q_{e}}\in (\mathfrak{g}\cdot q_{e})^{\perp }$ such that on
the open $G$-invariant neighborhood
$G\cdot \operatorname{Exp}_{q_{e}}(B)$ of $G\cdot q_{e}$, we
have $(G_{q_{e}})\preceq (G_{q})$ for every $q\in G\cdot
\operatorname{Exp}_{q_{e}}(B)$. Moreover
$G$ acts freely on $G\cdot \operatorname{Exp}_ {q_{e}}\left(B\cap
(T_{q_{e}}Q)_{\{e\}}\right)$. It is easy to see that $B\times
\mathfrak{g}^{\ast }$ can be identified with a slice at
$(q_{e},0)$ with respect to the diagonal action of $G$ on $
(G\cdot \operatorname{Exp}_{q_{e}}(B))\times \mathfrak{g}^{\ast
}$.
The strategy to prove the existence of a bifurcating branch of relative
equilibria with no symmetry from the set of relative equilibria
$\mathfrak{t}
\cdot q_{e}$ is the following. Note that we do not know a priori
which relative equilibrium in $\mathfrak{t} \cdot q_{e}$ will
bifurcate. We search for a local bifurcating branch of relative
equilibria in the following manner. Take a vector $v_{q_{e}}\in
B\cap (T_{q_{e}} Q)_{\{e\}}$ and note that
$\operatorname{Exp}_{q_e} (v_{q_{e}})
\in Q$ is a point with no symmetry, that is,
$G_{\operatorname{Exp}_{q_e} (v_{q_{e}})} = \{e\}$. Then $\tau
v_{q_{e}} \in B\cap (T_{q_{e}} Q)_{\{e\}}$, for $\tau\in
I $, where $I $ is an open interval containing $[0,1]$, and
$\operatorname{Exp}_{q_e} (\tau v_{q_{e}})$ is a smooth path
connecting $q_e $, the base point of the relative equilibrium in
$\mathfrak{t}\cdot q_e $ containing the branch of bifurcating
relative equilibria, to
$\operatorname{Exp}_{q_e} (v_{q_{e}})
\in Q$. In addition, we shall impose that the entire path
$\operatorname{Exp}_{q_e} (\tau v_{q_{e}})$ be formed by
base points of relative equilibria. We still need the vector part
of these relative equilibria which we postulate to be of the form
$\zeta(\tau)_Q(\operatorname{Exp}_{q_e} (\tau v_{q_{e}}) )$,
where $\zeta(\tau) \in \mathfrak{g}$ is a smooth path of Lie
algebra elements with $\zeta(0) \in \mathfrak{t}$. Since
$\operatorname{Exp}_{q_e} (\tau v_{q_{e}})$ has no symmetry for
$\tau> 0 $, the locked inertia tensor is invertible at these
points and the path
$\zeta(\tau) $ will be of the form
\[
\zeta(\tau) = \mathbb {I}(\operatorname{Exp}_{q_e} (\tau
v_{q_{e}}))^{-1} (\beta(\tau)),
\]
where $\beta(\tau) $ is a smooth path in $\mathfrak{g}^\ast$ with
$\beta(0) \in \mathbb {I}(q_e) \mathfrak{t}$. Now we shall use
the characterization of relative equilibria involving the amended
potential to require that the path
$\left(\operatorname{Exp}_{q_e}(\tau v_{q_e}), \beta(\tau) \right)
\in (G\cdot \operatorname{Exp}_{q_{e}}(B)) \times
\mathfrak{g}^\ast $ be such that
$f^{-1}(\left(\operatorname{Exp}_{q_e}(\tau v_{q_e}), \beta(\tau)
\right)$ are all relative equilibria.
The amended potential criterion is  applicable along the path
$\operatorname{Exp}_{q_e}(\tau v_{q_e})$ for $\tau >0 $, because
these points have no symmetry. As we shall see below, we shall
look for $\beta(\tau)$ of a certain form and
then the characterization of relative equilibria via  the amended
potential will impose conditions on both $\beta(\tau) $ and
$v_{q_e}$. We begin by specifying the form of $\beta(\tau)$.
\subsection{Splittings}
We shall need below certain direct sum decompositions of
$\mathfrak{g}$ and $\mathfrak{g}^\ast$.  The compactness of
$G$ implies that $\mathfrak{g}$ has an invariant inner
product and that $\mathfrak{g}=\mathfrak{t}\oplus
[\mathfrak{g},\mathfrak{t}]$ is an orthogonal direct sum. Let
$\mathfrak{k}_{1}\subset\mathfrak{t}$ be the orthogonal complement
to $\mathfrak{k}_{0}:=\mathfrak{g}_{q_{e}}$ in
$\mathfrak{t}$. Denoting
$\mathfrak{k}_{2}:=[\mathfrak{g},\mathfrak{t}]$ we obtain the
orthogonal direct sum $\mathfrak{g}=\mathfrak{k}_{0}\oplus
\mathfrak{k} _{1}\oplus \mathfrak{k}_{2}$.
For the dual of the Lie algebra, let
$\mathfrak{m}_{i}:=(\mathfrak{k} _{j}\oplus
\mathfrak{k}_{k})^{\circ }$ where $(i,j,k)$ is a cyclic
permutation of $(0,1,2)$. Then $\mathfrak{g}^{\ast }=\mathfrak{m}
_{0}\oplus \mathfrak{m}_{1}\oplus \mathfrak{m}_{2}$ is also an
orthogonal direct sum relative to the inner product on
$\mathfrak{g}^\ast$ naturally induced by the invariant inner
product on $\mathfrak{g}$.
\begin{lemma}
The subspaces defined by the above splittings have the following
properties:
\begin{enumerate}
\item[\textbf{(i)}] $\mathfrak{k}_{0}$, $\mathfrak{k}_{1}$,
$\mathfrak{k}_{2}$ are $G_{q_{e}}$-invariant and $G_{q_{e}}$ acts
trivially on $\mathfrak{k}_{0}$ and $\mathfrak{k}_{1}$;
\item[\textbf{(ii)}] $\mathfrak{m}_{0}$, $\mathfrak{m}_{1}$,
$\mathfrak{m}_{2}$ are $G_{q_{e}}$-invariant and $G_{q_{e}}$ acts
trivially on $\mathfrak{m}_{0}$ and $\mathfrak{m}_{1}$.
\end{enumerate}
\end{lemma}
\begin{proof}
\textbf{(i)} Because $G_{q_{e}}$ is a subgroup of $ \mathbb{T}$
it is obvious that $G_{q_{e}}$ acts trivially on
$\mathfrak{t}=\mathfrak{k}_{0}\oplus \mathfrak{k}_{1}$ and hence
on each summand. To prove the
$G_{q_{e}}$-invariance of
$\mathfrak{k}_{2}=[\mathfrak{g},\mathfrak{t}]
$, we use the fact that
$\operatorname{Ad}_{g}[\xi _{1},\xi
_{2}]=[\operatorname{Ad}_{g}\xi _{1},\operatorname{Ad}_{g}\xi
_{2}]$, for any
$\xi _{1},\xi _{2}\in
\mathfrak{g}$ and $g\in G$. Indeed, if  $\xi _{1}\in
\mathfrak{g}$,
$\xi _{2}\in \mathfrak{t}$, $g\in G_{q_{e}}$ we get
$\operatorname{Ad}_{g}[\xi _{1},\xi _{2}]\in
\lbrack \mathfrak{g},\mathfrak{t}]=\mathfrak{k}_{2}$.
\textbf{(ii)} For $g\in G_{q_{e}}$, $\mu \in \mathfrak{m}_{0}$ we
have to prove that $ \operatorname{Ad}_{g}^{\ast }\mu \in
\mathfrak{m}_{0}$. Indeed, if $\xi =\xi _{1}+\xi _{2}\in
\mathfrak{k}_{1}\oplus \mathfrak{k}_{2}$, we have
\begin{align*}
\langle \operatorname{Ad}_{g}^{\ast }\mu ,\xi \rangle
&=\langle \operatorname{Ad}_{g}^{\ast
}\mu ,\xi _{1}+\xi _{2}\rangle
=\langle \mu, \operatorname{Ad}_{g}(\xi _{1}+\xi _{2})\rangle \\
&=\langle \mu ,\xi _{1}+\operatorname{Ad}_{g}\xi _{2}\rangle
=0
\end{align*}
since  $G_{q_{e}}$ acts trivially on $\mathfrak{k}_{1}$,
$\mathfrak{k}_{2}$ is $G_{q_{e}}$--invariant and
$\mathfrak{m}_{0}=(\mathfrak{k}_{1}\oplus
\mathfrak{k}_{2})^{\circ }$.
The same type of proof holds for $\mathfrak{m}_{1}$and
$\mathfrak{m}_{2}$.
For $g\in G_{q_{e}}$, $\mu \in \mathfrak{m}_{0}$ we have to prove that $
\operatorname{Ad}_{g}^{\ast }\mu =\mu $.
Let $\xi =\xi _{0}+\xi _{1}+\xi _{2}\in \mathfrak{g}$, with
$\xi_i \in \mathfrak{k}_i $, $i = 0,1,2 $.
We have
\begin{align*}
\langle \operatorname{Ad}_{g}^{\ast }\mu -\mu ,\xi
\rangle
&=\langle \operatorname{Ad}_{g}^{\ast }\mu ,\xi _{0}+\xi _{1}+\xi
_{2}\rangle -\langle \mu ,\xi _{0}+\xi _{1}+\xi _{2}\rangle \\
&=\langle \mu ,\operatorname{Ad}_{g}(\xi _{0}+\xi _{1}+\xi
_{2})\rangle - \langle \mu ,\xi _{0}+\xi _{1}+\xi _{2}\rangle \\
&=\langle \mu ,\xi _{0}+\xi _{1}+\operatorname{Ad}_{g}\xi
_{2}\rangle - \langle \mu ,\xi _{0}\rangle
=\langle \mu ,\xi_{1}+\operatorname{Ad}_{g}\xi _{2} \rangle
=0
\end{align*}
because $G_{q_{e}}$ acts trivially on $\mathfrak{k}_{0}\oplus \mathfrak{k}
_{1}$, $\mathfrak{k}_{2}$ is $G_{q_{e}}$--invariant, and
$\mathfrak{m}_{0}=(
\mathfrak{k}_{1}\oplus \mathfrak{k}_{2})^{\circ }$.
The same type of proof holds for $\mathfrak{m}_{1}$.
\end{proof}
Recall from \S \ref{Simple mechanical systems with symmetry} that
$\ker \mathbb {I}(q_e) =
\mathfrak{g}_{q_e} = \mathfrak{k}_0 $. In particular, $\mathbb
{I}(q_e) \mathfrak{k}_0 = \{0\}$. The value of $\mathbb {I}(q_e)$
on the other summands in the decomposition $\mathfrak{g}=
\mathfrak{k}_0 \oplus \mathfrak{k}_1 \oplus \mathfrak{k}_2 $ is
given by the following lemma.
\begin{lemma}
\label{moment of  inertia isomorphism}
For $i\in \{1,2\}$ we have that $\mathfrak{m}_{i}=\mathbb{I}(q_{e})
\mathfrak{k}_{i}$.
\end{lemma}
\begin{proof}
Let $\kappa_{i}\in\mathfrak{k}_{i}$ with $i\in \{0,1,2\}$ be
arbitrary. Then
\begin{equation*}
\langle \mathbb{I}(q_{e})\kappa_{1},
\kappa_{0}+\kappa_{2}\rangle
=\langle \mathbb{I}(q_{e})\kappa_{1},\kappa_{0}\rangle
+\langle\mathbb{I}(q_{e})\kappa_{1},\kappa_{2}\rangle
=\langle \mathbb{I}(q_{e})\kappa_{0},\kappa_{1}\rangle
+\langle\mathbb{I}(q_{e})\kappa_{1},\kappa_{2}\rangle
=0
\end{equation*}
as $\ker \mathbb{I}(q_{e})=\mathfrak{k}_{0}$ and, by Proposition
\ref{montaldi} (ii),
$\mathbb{I}(q_{e})\mathfrak{t}\subset \mathfrak{k}_{2}^{\circ}$.
This proves that $\mathbb{I}(q_{e})\mathfrak{k}_{1}\subset
\mathfrak{m}_{1}$. Counting dimensions we have that
$\dim\mathbb{I}(q_{e})\mathfrak{k}_{1}= \dim \mathfrak{k}_1 -
\dim \ker \left(\mathbb {I}(q_e)|_{\mathfrak{k}_1} \right) =
\dim \mathfrak{g}-\dim \mathfrak{k}_{0}-\dim \mathfrak{k}_{2}
=\dim \mathfrak{m}_{1}$,
since $\ker \left(\mathbb {I}(q_e)|_{\mathfrak{k}_1} \right) =
\{0\} $. This proves that
$\mathfrak{m}_{1}=\mathbb{I}(q_{e})
\mathfrak{k}_{1}$. In an analogous way we prove the equality for
$i=2$.
\end{proof}
In the next paragraph we shall need the direct sum decomposition
$\mathfrak{g}^\ast =  \mathfrak{m}_1 \oplus \mathfrak{m}$, where
$\mathfrak{m}_1 = \mathbb {I}(q_e) \mathfrak{t}$ and
$\mathfrak{m} := \mathfrak{m}_0 \oplus \mathfrak{m}_2 $.
Let $\Pi_1: \mathfrak{g}^\ast\rightarrow \mathbb {I}(q_e)
\mathfrak{t}$ be the projection along $\mathfrak{m} $. Similarly,
denote
$\mathfrak{k}: =
\mathfrak{k}_1
\oplus\mathfrak{k}_2$, and write $\mathfrak{g} =
\mathfrak{g}_{q_e} \oplus\mathfrak{k}$. Thus there is another
decomposition of $\mathfrak{g}^\ast$, namely, $\mathfrak{g}^\ast
= \mathfrak{g}_{q_e} ^\circ \oplus \mathfrak{k}^\circ $. However,
for any $\zeta \in \mathfrak{g}_{q_e} $ and any $\xi\in
\mathfrak{g}$,  we have
$\langle \mathbb {I}(q_e) \xi, \zeta \rangle = \langle\!\langle
\xi_Q(q_e), \zeta_Q(q_e) \rangle\!\rangle = 0 $ since
$\zeta_Q(q_e) = 0 $, which shows that $\mathbb {I}(q_e)
\mathfrak{g} \subset \mathfrak{g}_{q_e}^\circ $. Since $ \ker
\mathbb {I}(q_e) = \mathfrak{g}_{q_e} $, it follows that $\dim
\mathbb {I}(q_e) \mathfrak{g} = \dim \mathfrak{g} - \dim \ker
\mathbb {I}(q_e) = \dim \mathfrak{g} - \dim \mathfrak{g}_{q_e} =
\dim \mathfrak{g}_{q_e}^\circ $, which shows that
$\mathfrak{g}_{q_e}^\circ = \mathbb {I}(q_e)\mathfrak{g}$. Thus
we also have the direct sum decomposition $\mathfrak{g}^\ast =
\mathbb {I}(q_e)\mathfrak{g} \oplus \mathfrak{k}^\circ $. Note
that $\mathbb {I}(q_e) \mathfrak{g} = \mathfrak{m}_1 \oplus
\mathfrak{m}_2 $, by Lemma \ref{moment of  inertia isomorphism}
and that $\mathfrak{m}_0 = \mathfrak{k}^\circ $. Summarizing we
have:
\[
\mathfrak{g}^\ast = \mathfrak{m}_0 \oplus \mathfrak{m}_1 \oplus
\mathfrak{m}_2 =  \mathfrak{k}^\circ \oplus \mathbb {I}(q_e)
\mathfrak{g}, \quad \text{where} \quad \mathbb {I}(q_e) \mathfrak{g} = \mathfrak{m}_1 \oplus
\mathfrak{m}_2 \quad \text{and} \quad \mathfrak{m}_0 =
\mathfrak{k}^\circ.
\]
\subsection{The rescaled equation}
Recall that $B\subset (\mathfrak{g}\cdot q_{e})^{\perp }$ is a
$G_{q_{e}}$-invariant open neighborhood of
$0_{q_{e}}\in (\mathfrak{g}\cdot q_{e})^{\perp }$ such that on
the open $G$-invariant neighborhood
$G\cdot \operatorname{Exp}_{q_{e}}(B)$ of $G\cdot q_{e}$, we
have $(G_{q_{e}})\preceq (G_{q})$ for every $q\in G\cdot
\operatorname{Exp}_{q_{e}}(B)$.
Consider the following rescaling:
\begin{equation*}
v_{q_{e}}\in B\cap (T_{q_{e}}Q)_{\{e\}}\mapsto \tau v_{q_{e}}\in
B\cap (T_{q_{e}}Q)_{\{e\}}
\end{equation*}
\begin{equation*}
\mu \in \mathfrak{g}^{\ast }\mapsto \beta (\tau ,\mu )\in \mathfrak{g}^{\ast }
\end{equation*}
where, $\tau \in I$, $I$ is an open interval containing
$[0,1]$, and
$\beta :I\times \mathfrak{g}^{\ast }\rightarrow
\mathfrak{g}^{\ast }$ is chosen such that $\beta (0,\mu
)=\Pi _{1}\mu $. So, for $(v_{q_{e}},\mu )$ fixed, $(\tau
v_{q_{e}},\beta (\tau ,\mu ))$ converges to $(0_{q_{e}},\Pi
_{1}\mu )$ as $\tau \rightarrow 0$. Define
\begin{equation*}
\beta (\tau ,\mu ):=\Pi _{1}\mu +\tau \beta '(\mu
)+\tau^{2} \beta ''(\mu )
\end{equation*}
for some arbitrary smooth functions $\beta ',\beta '':
\mathfrak{g}^{\ast }\rightarrow \mathfrak{g}^{\ast }$. Since
$\mathbb{I}$ is invertible only for points with no symmetry, we
want to find conditions on
 $\beta '$, $\beta ''$ such that the expression
\begin{equation}
\label{32}
\mathbb{I}(\operatorname{Exp} _{q_{e}}(\tau v_{q_{e}}))^{-1}\beta
(\tau ,\mu )
\end{equation}
extends to a smooth function in a neighborhood of $\tau =0$. Note
that $v_{q_{e}}$ is different from $0_{q_{e}}$ since
$G_{v_{q_{e}}} = \{e\}$ by construction and $G_{0_{q_{e}}} =
G_{q_e}
\neq \{e\}$. Define
\begin{equation*}
\Phi :I\times \left(B\cap (T_{q_{e}}Q)_{\{e\}}\right)\times
\mathfrak{g}^{\ast }\times
\mathfrak{g}_{q_{e}}\times \mathfrak{k}\rightarrow \mathfrak{g}^{\ast }
\end{equation*}
\begin{equation}
\label{definition of Phi}
\Phi (\tau ,v_{q_{e}},\mu ,\xi ,\eta )
:=\mathbb{I}(\operatorname{Exp}_{q_{e}}(\tau
v_{q_{e}})) (\xi +\eta )-\beta (\tau ,\mu ).
\end{equation}
Now we search for the velocity $\xi + \eta $ of relative
equilibria among the solutions of
$\Phi (\tau ,v_{q_{e}},\mu ,\xi ,\eta )=0$. We shall prove below
that $\xi $ and $\eta$ are smooth functions of
$\tau $, $v_{q_{e}}$, $\mu$, even at $\tau= 0$.
Then \eqref{32} shows that $\xi+
\eta $ is a smooth function of $\tau $, $v_{q_{e}}$, $\mu$, for
$\tau$ in a small neighborhood of zero.
\subsection{The Lyapunov-Schmidt procedure}
To solve $\Phi =0$ we apply the standard Lyapunov-Schmidt
method. This equation has  a unique solution for
$\tau
\neq 0$, because $\tau v_{q_{e}}\in B\cap (T_{q_{e}}Q)_{\{e\}}$ so
$\mathbb{I}(\operatorname{Exp}_{q_{e}}(\tau v_{q_{e}}))$ is
invertible. It remains to prove that the equation has a solution
when $\tau =0$.
Denote by $D_{\mathfrak{g}_{q_{e}}\times \mathfrak{k}}$ the
Fr\'echet derivative relative to the last two factors
$\mathfrak{g}_{q_{e}}\times \mathfrak{k}$  in the definition of
$\Phi$. We have
\begin{equation*}
\ker D_{\mathfrak{g}_{q_{e}}\times \mathfrak{k}}\Phi (0,v_{q_{e}},\mu ,\xi
,\eta )=\ker \mathbb{I}(q_{e})=\mathfrak{g}_{q_{e}}.
\end{equation*}
We will solve the equation $\Phi =0$ in two steps. For this, let
\begin{equation*}
\Pi :\mathfrak{g}^{\ast }\rightarrow \mathbb{I}(q_{e})\mathfrak{g}
\end{equation*}
be the projection induced by the splitting $\mathfrak{g}^{\ast }=\mathbb{I}
(q_{e})\mathfrak{g}\oplus \mathfrak{k}^{\circ }$.
\textbf{Step1}. Solve $\Pi \circ \Phi =0$ for $\eta $ in terms of
$\tau $, $v_{q_{e}}$, $\mu $, $\xi $. For this, let
\begin{align*}
\widehat{\mathbb{I}}(\operatorname{Exp} _{q_{e}}(\tau v_{q_{e}}))
&:= (\Pi\circ \mathbb{I})(\operatorname{Exp}_{q_{e}}
(\tau v_{q_{e}}))|_{\mathfrak{k}}: \mathfrak{k}\rightarrow
\mathbb {I}(q_e) \mathfrak{g}\\
\overset{\thicksim}{\mathbb{I}}(\operatorname{Exp} _{q_{e}}
(\tau v_{q_{e}}))
&:=(\Pi \circ \mathbb{I})(\operatorname{Exp}_{q_{e}}(\tau
v_{q_{e}}))|_{\mathfrak{g}_{q_{e}}}: \mathfrak{g}_{q_e}
\rightarrow \mathbb {I}(q_e) \mathfrak{g}
\end{align*}
where $\widehat{\mathbb{I}}(\operatorname{Exp} _{q_{e}}(\tau
v_{q_{e}})) $ is an isomorphism even when $\tau =0$. Then we
obtain
\begin{equation}
\label{pi composed with Phi}
(\Pi \circ \Phi )(0,v_{q_{e}},\mu ,\xi ,\eta )
= \Pi [\mathbb{I} (q_{e})(\xi +\eta )-\beta (0,\mu)]
=\widehat{\mathbb{I}}(q_{e}) \eta -\Pi _{1}\mu .
\end{equation}
Denoting $\eta _{\mu }:=\widehat{\mathbb{I}}
(q_{e})^{-1}(\Pi _{1}\mu )$, we have $(\Pi \circ \Phi)
(0,v_{q_{e}},\mu ,\xi ,\eta _{\mu })\equiv 0$. 
Denoting by $D_\eta$ the partial Fr\'echet derivative relative
to the variable $\eta
\in \mathfrak{k}$ we get at any given point $(0, v_{q_e}^0,
\mu^0, \xi^0, \eta^0 )$
\begin{equation}
\label{isomorphism condition for ift}
D_\eta (\Pi \circ \Phi )(0, v_{q_e}^0, \mu^0, \xi^0, \eta^0)
=\widehat{\mathbb{I}}(q_{e})
\end{equation}
which is invertible. Thus the implicit function theorem gives a
unique smooth function
$\eta (\tau ,v_{q_{e}},\mu ,\xi )$ such that $\eta
(0, v_{q_e}^0,
\mu^0, \xi^0)=\eta^0$ and
\begin{equation}
\label{definition of eta}
(\Pi \circ \Phi) (\tau ,v_{q_{e}},\mu ,\xi ,\eta (\tau
,v_{q_{e}},\mu ,\xi ))\equiv 0.
\end{equation}
The function $\eta $ is defined in some open set in $I\times
\left(B\cap (T_{q_{e}}Q)_{\{e\}}\right)\times
\mathfrak{g}^{\ast}\times
\mathfrak{g}_{q_{e}}$ containing $(0, v_{q_e}^0, \mu^0, \xi^0)
\in \{0\}\times
\left(B\cap (T_{q_{e}}Q)_{\{e\}}\right)\times
\mathfrak{g}^{\ast }\times \mathfrak{g}_{q_{e}}$.
If we now choose $\eta^0 = \eta_{\mu^0} = \widehat{\mathbb
{I}}(q_e) ^{-1}( \Pi_1 \mu ^0)$, then uniqueness of the
solution of the implicit function theorem implies that $\eta(0,
v_{q_e}, \mu, \xi) = \eta_\mu$ in the
neighborhood of $(0, v_{q_e}^0, \mu^0, \xi^0)$.
\medskip
Later we will need the following result.
\begin{proposition}
\label{belongs}
We have $\eta _{\mu }:=\widehat{\mathbb{I}}
(q_{e})^{-1}(\Pi _{1}\mu ) \in \mathfrak{k}_1
\subset \mathfrak{t}$.
\end{proposition}
\begin{proof}
Since we can write $\mathfrak{t}=\ker \mathbb{I}(q_{e})\oplus
\mathfrak{k}_1$ we obtain
\begin{equation*}
\widehat{\mathbb{I}}(q_{e})\mathfrak{k}_1
=(\Pi \circ \mathbb{I}(q_{e}))
\mathfrak{k}_1
=\mathbb{I} (q_{e})
\mathfrak{k}_1=\mathbb{I}(q_{e})(\mathfrak{t})
=\operatorname{Im}\Pi _{1}.
\end{equation*}
Now, because $\widehat
{\mathbb{I}}(q_{e})$ is an isomorphism, it follows that
$\widehat{\mathbb{I}}
(q_{e})^{-1}(\Pi _{1}\mu ) \in \mathfrak{k}_1$.
\end{proof}
\textbf{Step2}. Now we solve the equation $(Id-\Pi)\circ\Phi=0$.
For this, let
\begin{equation*}
\varphi :I\times \left(B\cap (T_{q_{e}}Q)_{\{e\}}\right)\times
\mathfrak{g}^{\ast }\times
\mathfrak{g}_{q_{e}}\rightarrow \mathfrak{k}^{\circ }
\end{equation*}
\begin{equation}
\label{definition of phi}
\varphi (\tau ,v_{q_{e}},\mu ,\xi )
:=(Id-\Pi ) \Phi (\tau ,v_{q_{e}},\mu ,\xi ,\eta (\tau
,v_{q_{e}},\mu ,\xi )).
\end{equation}
In particular, $\varphi (0,v_{q_{e}},\mu ,\xi )=(Id-\Pi)
(\mathbb{I} (q_{e})(\xi +\eta _{\mu })-\Pi _{1}\mu )$. Since
$\operatorname{Im}\mathbb{I} (q_{e}) = \operatorname{Im}\Pi$ and
$\operatorname{Im}\Pi _{1} = \mathbb {I}(q_e) \mathfrak{t}
\subset \mathbb {I}(q_e) \mathfrak{g}$, it follows that $\varphi
(0,v_{q_{e}},\mu ,\xi )\equiv 0$. We shall solve for $\xi
\in\mathfrak{g}_{q_e}$, in the
neighborhood of $(0, v_{q_e}^0, \mu^0, \xi^0)$ found in Step 1,
the equation $\varphi (\tau ,v_{q_{e}},\mu ,\xi ) = 0 $. To do
this, we shall need information about the higher derivatives of 
$\varphi $ with respect to $\tau $, evaluated at $\tau =0$.
\begin{lemma}
Let $\xi $, $\eta \in \mathfrak{g}$ and $q\in Q$. Suppose that
$\mathbf{d}V_{\eta }(q)=0$, where $V_{\eta}$ is the augmented
potential and suppose that both $\xi $ and $[\xi ,\eta ]$ belong
to $\mathfrak{g}_{q}$. Then
$\mathbf{d}\langle \mathbb{I} (\cdot )\xi ,\eta \rangle (q)=0$.
\end{lemma}
\begin{proof}
Since $\mathbf{d}V_{\eta}(q)=0$, $\eta _{Q}(q)$ is a relative
equilibrium by Proposition \ref{augmented potential criterion},
that is,
$X_{H}(\alpha _{q})=\eta_{T^{\ast }Q}(\alpha _{q})$, where $\alpha
_{q}=\mathbb{F}L(\eta _{Q}(q))$. Now suppose that both $\xi $,
$[\xi ,\eta ]\in \mathfrak{g}_{q}$. Then
\begin{equation*}
\xi _{T^{\ast }Q}(\alpha _{q})=\left. \frac{d}{dt}\right|
_{t=0}\mathbb{F}L(\exp(t\xi )\cdot \eta
_{Q}(q))=\mathbb{F}L([\xi ,\eta ]_{Q}(q))=0,
\end{equation*}
where we have used that $g\cdot \eta
_{Q}(q)=(\operatorname{Ad}_{g}\eta )_{Q}(g\cdot q)$. It follows
that
$(\eta +\xi )_{T^{\ast }Q}(\alpha _{q})=X_{H}(\alpha _{q})$ and
hence, again by Proposition \ref{augmented potential criterion},
that $0=\mathbf{d}V_{\eta +\xi }(q)=\mathbf{d}V_{\eta}(q)
-\mathbf{d}\langle\mathbb{I}(\cdot )\eta ,\xi
\rangle(q)-\frac{1}{2}\mathbf{d}\|
\xi _{Q}(\cdot )\| ^{2}(q)$. However, $\mathbf{d}\| \xi
_{Q}(\cdot )\| ^{2}(q) = 0 $ since $\xi \in \mathfrak{g}_q $, as
an easy coordinate computation shows. Since $\mathbf{d}V_\eta(q) =
0 $ by hypothesis, we have $\mathbf{d}\langle\mathbb{I}(\cdot
)\eta ,\xi \rangle(q) = 0 $. Symmetry of $\mathbb {I}(q)$ proves
the result.
\end{proof}
Let now $\xi \in \mathfrak{g}_{q_{e}}$ and $\eta \in
\mathfrak{t}$. Since $\mathfrak{g}_{q_e} \subset \mathfrak{t}$,
we have $[ \xi, \eta ] = 0 \in \mathfrak{g}_{q_e}$. In
addition, hypothesis \textbf{(H)} and Proposition
\ref{augmented potential criterion}, guarantee that
$\mathbf{d}V_\xi(q_e) = 0 $ which shows that all hypotheses of
the previous lemma are satisfied. Therefore,
\begin{equation}
\label{differential of i}
\mathbf{d}\langle\mathbb{I}(\cdot )\xi,\eta
\rangle(q_{e})=0 \quad  \text{for} \quad \xi \in
\mathfrak{g}_{q_{e}}, \; \eta \in
\mathfrak{t}.
\end{equation}
\subsection{The bifurcation equation}
Now we can proceed with the study of equation $
\varphi= (Id-\Pi)\circ\Phi=0$.
We have
\begin{align}
\label{first tau derivative of phi}
\frac{\partial \varphi}{\partial \tau}(\tau, v_{q_e}, \mu, \xi)
= (Id - \Pi) &\left[ T_{\tau v_{q_e}}(\mathbb {I} \circ
\operatorname{Exp}_{q_e})(v_{q_e})(\xi + \eta(\tau, v_{q_e}, \mu,
\xi))
+ \mathbb {I}(\operatorname{Exp}_{q_e}(\tau v_{q_e}))
\frac{\partial \eta}{\partial \tau}(\tau, v_{q_e}, \mu, \xi)
\right. \nonumber \\
&\qquad \left. - \frac{\partial \beta}{\partial \tau}(\tau, \mu)
\right].
\end{align}
\begin{proposition}
$\frac{\partial }{\partial \tau }\varphi (0,v_{q_{e}},\mu ,\xi )\equiv
-(Id-\Pi )\beta '(\mu )$.
\end{proposition}
\begin{proof}
Formula  \eqref{first tau derivative of phi} gives for $\tau= 0 $
\begin{equation*}
\frac{\partial \varphi }{\partial \tau }(0,v_{q_{e}},\mu ,\xi )
=(Id-\Pi ) \left[\big(T_{q_e}
\mathbb{I}(v_{q_{e}})\big)(\xi +\eta _{\mu })
+\mathbb{I}(q_{e})\frac{\partial \eta }{\partial
\tau }(0,v_{q_{e}},\mu ,\xi )-\frac{\partial \beta }{\partial \tau
}(0,\mu )\right].
\end{equation*}
Now, because $\operatorname{Im}\mathbb{I}(q_{e})=\operatorname{Im}\Pi $
we obtain $(Id-\Pi )
\circ \mathbb{I}(q_{e})=0$ and hence the second summand vanishes.
From
\eqref{differential of i} we have that
$(T_{q_e}\mathbb{I}
(v_{q_{e}}))(\mathfrak{t})\subset \mathfrak{g}_{q_{e}}^{\circ
}=\operatorname{Im}\Pi $. Using Proposition \ref{belongs} and
since $\xi \in \mathfrak{g}_{q_{e}}\subset
\mathfrak{t}$, we obtain that $\xi +\eta_{\mu }\in \mathfrak{t}$. Therefore
$(Id-\Pi )[(T_{q_e}\mathbb{I}(v_{q_{e}}))(\xi +\eta
_{\mu })]=0$.  Since $\frac{\partial \beta }{\partial \tau }(0,\mu
)=\beta ^{\prime }(\mu )$,  we obtain the desired equality.
\end{proof}
Let us impose the additional condition $\beta '(\mu)
\subset\operatorname{Im}\Pi $. Then it follows that
\begin{equation*}
\varphi (\tau ,v_{q_{e}},\mu ,\xi )=\tau ^{2}\psi (\tau ,v_{q_{e}},\mu ,\xi).
\end{equation*}
for some smooth function $\psi $ where
\begin{equation*}
\psi (0,v_{q_{e}},\mu ,\xi )=\frac{1}{2}\frac{\partial ^{2}\varphi }{
\partial \tau ^{2}}(0,v_{q_{e}},\mu ,\xi )
\end{equation*}
We begin by solving the equation
\begin{equation*}
\psi (0,v_{q_{e}},\mu ,\xi )=0
\end{equation*}
for $\xi $ as a function of $v_{q_{e}}$ and $\mu $. Equivalently,
we have to  solve
\begin{equation*}
\frac{1}{2}\frac{\partial ^{2}\varphi }{\partial \tau ^{2}}(0,v_{q_{e}},\mu
,\xi )=0.
\end{equation*}
To compute this second derivative of $\varphi$ we shall use
\eqref{first tau derivative of phi}. We begin by noting that
$\tau\in I \mapsto T_{\tau v_{q_e}}(\mathbb {I} \circ
\operatorname{Exp}_{q_e})(v_{q_e})$ is a smooth path in
$L(\mathfrak{g}, \mathfrak{g}^\ast)$ and so we can define the
linear operator from $\mathfrak{g}$ to $\mathfrak{g}^\ast$
by
\[
A_{v_{q_e}}: = \left.\frac{\partial}{\partial \tau}\right|_{\tau=
0}T_{\tau v_{q_e}}(\mathbb {I} \circ
\operatorname{Exp}_{q_e})(v_{q_e}) \in L(\mathfrak{g},
\mathfrak{g}^\ast).
\]
With this notation, formulas \eqref{first tau derivative of phi},
\eqref{definition of Phi}, \eqref{definition of phi}, and
Proposition \ref{belongs} yield
\begin{align}
\label{second tau derivative of phi at zero}
\frac{\partial ^{2}\varphi }{\partial \tau^{2}} (0,v_{q_{e}},\mu,
\xi )
&=(Id-\Pi)\left[ A_{v_{q_e}}(\xi +\eta _{\mu})
+2T_{q_e}\mathbb{I}(v_{q_{e}})\frac{\partial \eta}{\partial
\tau }(0,v_{q_{e}},\mu ,\xi ) \right.\\
& \qquad \qquad \qquad \left. + \mathbb{I}(q_{e})\frac{\partial
^{2}\eta }{\partial
\tau ^{2}} (0,v_{q_{e}},\mu ,\xi )-2\beta ''(\mu
)\right] \nonumber \\
&= (Id-\Pi)\left[ A_{v_{q_e}}(\xi +\eta _{\mu})
+2T_{q_e}\mathbb{I}(v_{q_{e}})\frac{\partial \eta}{\partial
\tau }(0,v_{q_{e}},\mu ,\xi )  - 2\beta ''(\mu) \right]
\nonumber
\end{align}
since $(Id - \Pi)\mathbb{I}(q_{e})\frac{\partial
^{2}\eta }{\partial \tau ^{2}} (0,v_{q_{e}},\mu ,\xi) = 0$.
Let $\{\xi _{1},...,\xi _{p}\}$ be a basis of
$\mathfrak{g}_{q_{e}}$. Since $\partial^2 \varphi ( \tau, v_{q_e},
\mu, \xi)/\partial \tau ^2 \in \mathfrak{k}^\circ $ and
$\mathfrak{g} = \mathfrak{g}_{q_e} \oplus \mathfrak{k}$ the
equation $\partial ^2 \varphi(0, v_{q_e}, \mu, \xi)/\partial
\tau^2 = 0 $ is equivalent to the following system of $p $
equations
\[
\left\langle\frac{\partial ^{2}\varphi }{\partial \tau ^{2}}(0,v_{q_{e}},\mu
,\xi ) , \xi_b \right \rangle = 0, \quad \text{for~all} \quad b =
1, \dots, p,
\]
which, by \eqref{second tau derivative of phi at zero}, is
\[
\left\langle (Id-\Pi)\left[ A_{v_{q_e}}(\xi +\eta _{\mu})
+2T_{q_e}\mathbb{I}(v_{q_{e}})\frac{\partial \eta}{\partial
\tau }(0,v_{q_{e}},\mu ,\xi )  - 2\beta ''(\mu
)\right] ,\xi_b \right\rangle = 0, \quad \text{for~all}
\quad b = 1, \dots, p.
\]
We shall show that in this expression we can drop the projector
$ Id - \Pi $. Indeed, let $\alpha = \alpha_0 + \alpha_1 +
\alpha_2 \in \mathfrak{g}^\ast = \mathfrak{m}_0 \oplus
\mathfrak{m}_1 \oplus \mathfrak{m}_2 $, where $\alpha_i \in
\mathfrak{m}_i $, for $i = 0,1,2 $. Since $\Pi:
\mathfrak{g}^\ast \rightarrow \mathbb {I}(q_e) \mathfrak{g} =
\mathfrak{m}_1 \oplus \mathfrak{m}_2 $, we have
\[
\langle (Id - \Pi) \alpha, \xi_b \rangle = \langle \alpha, \xi_b
\rangle - \langle \alpha_1, \xi_b \rangle - \langle \alpha_2,
\xi_b \rangle = \langle \alpha, \xi_b \rangle
\]
because $\langle \alpha_1, \xi_b \rangle = 0$, since $\alpha_1 \in
\mathfrak{m}_1 = (\mathfrak{k}_0 \oplus \mathfrak{k}_2)^\circ $,
$\xi_b
\in \mathfrak{g}_{q_e} = \mathfrak{k}_0 $, and
$\langle \alpha_2, \xi_b \rangle = 0 $, since $\alpha_2 \in
\mathfrak{m}_2 = (\mathfrak{k}_0 \oplus \mathfrak{k}_1)^\circ $,
$\xi_b
\in \mathfrak{g}_{q_e} = \mathfrak{k}_0 $. The system to be
solved is hence
\begin{equation}
\label{linear system to be solved}
\left\langle  A_{v_{q_e}}(\xi +\eta _{\mu})
+2T_{q_e}\mathbb{I}(v_{q_{e}})\frac{\partial \eta}{\partial
\tau }(0,v_{q_{e}},\mu ,\xi )  - 2\beta ''(\mu
) ,\xi_b \right\rangle = 0, \quad \text{for~all}
\quad b = 1, \dots, p.
\end{equation}
In what follows we need the expression for $\frac{\partial \eta}
{\partial \tau }(0,v_{q_{e}},\mu ,\xi)$.  Differentiating
\eqref{definition of eta} relative to $\tau$ at zero and taking
into account \eqref{isomorphism condition for ift} and
\eqref{definition of Phi}, we get
\begin{align}
\label{first tau derivative of eta at zero}
&\frac{\partial \eta }{\partial \tau }(0,v_{q_{e}},\mu ,\xi)
= - \widehat{\mathbb {I}}(q_e)^{-1}\frac{\partial}{\partial
\tau}(\Pi\circ \Phi)(0, v_{q_e}, \mu, \xi, \eta_\mu)\\
&\qquad = - \widehat{\mathbb {I}}(q_e)^{-1} \Pi\left[
T_{q_e}\mathbb {I}(v_{q_e})(\xi+ \eta_\mu) - \beta'(\mu) \right]
\nonumber \\
&\qquad = - \left(
\widehat{\mathbb{I}}(q_e)^{-1} \circ T_{q_e}\overset{\thicksim}
{\mathbb {I}}(v_{q_e})\right) \xi
- \left(\widehat{\mathbb
{I}}(q_e)^{-1} \circ  T_{q_e}\widehat{\mathbb
{I}}(v_{q_e}) \circ
\widehat{\mathbb{I}} (q_{e})^{-1}\right)(\Pi _{1}\mu ) +
\widehat{\mathbb{I}} (q_{e})^{-1}(\beta'(\mu)) \nonumber
\end{align}
since $T_{q_e}\overset{\thicksim}{\mathbb {I}} = \Pi \circ
T_{q_e}\mathbb {I}|_{\mathfrak{g}_{q_e}} $ and
$T_{q_e}\widehat{\mathbb {I}} = \Pi \circ
T_{q_e}\mathbb {I}|_{\mathfrak{k}}$.
Expanding
$\xi$ in the basis
$\{\xi_1,
\dots, \xi_p\}$ as $\xi= \alpha^i \xi_i $ and taking into account
the above expression, the system \eqref{linear system to be
solved} is equivalent to the following system of linear
equations in the unknowns  $\alpha^1, \dots, \alpha^p$
\[
A_{a b} \alpha^a + B_b = 0, \quad a,b = 1, \dots , p,
\]
where
\begin{align}
A_{a b}&: = \left\langle A_{v_{q_e}} \xi_a, \xi_b \right\rangle -
2 \left\langle \left(T_{q_e}\mathbb {I}(v_{q_e}) \circ
\widehat{\mathbb {I}}(q_e)^{-1} \circ T_{q_e} \overset{\thicksim}
{\mathbb {I}}(v_{q_e})\right) \xi_a, \xi_b
\right\rangle
\label{A}\\
B_b&:=\left\langle \left(A_{v_{q_e}} \circ \widehat{\mathbb{I}}
(q_{e})^{-1}\circ \Pi_1\right) \mu, \xi_b \right\rangle
- 2 \left\langle \left(T_{q_e}\mathbb {I}(v_{q_e}) \circ
\widehat{\mathbb
{I}}(q_e)^{-1} \circ  T_{q_e}\widehat{\mathbb
{I}}(v_{q_e}) \circ
\widehat{\mathbb{I}} (q_{e})^{-1} \circ \Pi_1 \right)\mu, \xi_b
\right\rangle
\label{B} \\
&\qquad  +2 \left\langle \left(T_{q_e}\mathbb {I}(v_{q_e}) \circ
\widehat{\mathbb{I}}(q_e)^{-1}\right)\beta'(\mu), \xi_b
\right\rangle
- \left\langle \beta''(\mu), \xi_b \right\rangle.
\nonumber
\end{align}
Denote by $A: = [A_{a b}] $ the $p\times p $ matrix with entries
$A_{a b} $. Thus, if $v_{q_{e}}\notin
\mathcal{Z}_{\mu}=:\{v_{q_{e}}\in B\cap (T_{q_{e}}Q)_{\{e\}}\mid
\det A=0\}$ this linear system has a unique solution for
$\alpha^1, \dots , \alpha^p $, that is for $\xi$, as function of $
v_{q_{e}}$, $\mu $. we shall denote this solution by
$\xi_0(v_{q_e}, \mu)$.
\textit{Summarizing, if $v_{q_e} \notin
\mathcal{Z}_{\mu}$, then $\xi_0(v_{q_e}, \mu) $ is the unique
solution of the equation} 
\begin{equation}
\label{second step in LS}
\frac{\partial ^{2}\varphi }{\partial
\tau ^{2}}(0,v_{q_{e}},\mu ,\xi )=0.
\end{equation}
\begin{lemma}
\label{properties of Z} The set $\mathcal{Z}_{\mu}$ is closed and
$G_{q_{e}}$--invariant in $B\cap (T_{q_e}Q)_{\{e\}}$.
\end{lemma}
\begin{proof}
The set $\mathcal{Z}_{\mu}$ is obviously closed. Since
$\mathfrak{k}$ is $G_{q_{e}}$--invariant it follows that
$\mathfrak{k}^\circ$ is $G_{q_e}$--invariant. Formula
\eqref{equivariance of I} shows that
$\mathbb{I}(q_{e})\mathfrak{g}$ is also $G_{q_e}$--invariant.
Thus the direct sum $\mathbb{I}(q_{e}) \mathfrak{g}\oplus
\mathfrak{k}^{\circ }$ is a $G_{q_{e}}$--invariant decomposition
of $\mathfrak{g}^{\ast }$ and therefore $\Pi : \mathfrak{g}^\ast
\rightarrow \mathbb {I}(q_e) \mathfrak{g}$ is
$G_{q_{e}}$--equivariant. From the $G_{q_{e}}$--equivariance of
$\operatorname{Exp}_{q_{e}}$ and \eqref{equivariance of I}, it
follows that $\mathbb{I} (\operatorname{Exp} _{q_{e}}(h\cdot
v_{q_{e}}))=\operatorname{Ad}^\ast_{h^{-1}}\circ
\mathbb{I}(\operatorname{Exp} _{q_{e}}(v_{q_{e}})) \circ
\operatorname{Ad}_{h^{-1}} = \operatorname{Ad}^\ast_{h^{-1}}\circ
\mathbb{I}(\operatorname{Exp} _{q_{e}}(v_{q_{e}}))$ for any $h
\in G_{q_e}$ since $G_{q_e} \subset \mathbb{T}$ and is therefore
Abelian. Thus
\begin{align*}
\overset{\thicksim}{\mathbb{I}}(\operatorname{Exp}
_{q_{e}}(h\cdot v_{q_{e}}))
&= \Pi \circ \mathbb{I}(\operatorname{Exp}
_{q_{e}}(T_{q_{e}}\Psi _{h}\cdot v_{q_{e}}))|_{\mathfrak{g}_{q_e}}
= \Pi\circ \operatorname{Ad}^\ast_{h^{-1}}\circ
\mathbb{I}(\operatorname{Exp}
_{q_{e}}(v_{q_{e}}))|_{\mathfrak{g}_{q_e}}\\
&= \operatorname{Ad}^\ast_{h^{-1}} \circ
\Pi\circ \mathbb{I}(\operatorname{Exp}
_{q_{e}}(v_{q_{e}}))|_{\mathfrak{g}_{q_e}}
= \operatorname{Ad}^\ast_{h^{-1}} \circ
\overset{\thicksim}{
\mathbb{I}}(\operatorname{Exp} _{q_{e}}(v_{q_{e}}))
\end{align*}
for all $h\in G_{q_{e}}$ and $v_{q_{e}}\in B$. Replacing here
$v_{q_e} $ by $sv_{q_e} $ and taking the $s $--derivative at
zero, shows that
$T_{q_e}\overset{\thicksim}{\mathbb {I}}(h\cdot v_{q_e}) \xi =
\operatorname{Ad}^\ast_{h^{-1}}\left(
T_{q_e}\overset{\thicksim}{\mathbb {I}}(v_{q_e}) \xi
\right)$ for any $h\in G_{q_{e}}$ and $\xi \in
\mathfrak{g}_{q_{e}}$, that is,
$T_{q_e}\overset{\thicksim}{\mathbb {I}}(v_{q_e}) \xi $
is $G_{q_{e}}$--equivariant
as a function of $v_{q_{e}}$, for all $\xi \in
\mathfrak{g}_{q_{e}}$. Similarly $ T_{q_e} \mathbb
{I}(h\cdot v_{q_e}) = \operatorname{Ad}^\ast_{h^{-1}} \circ
T_{q_e} \mathbb{I}( v_{q_e}) \circ \operatorname{Ad}_{h^{-1}}$.
 From \eqref{equivariance of I} and the definition of
$\widehat{\mathbb {I}}(q_e)^{-1}$, it follows that
$\widehat{\mathbb {I}}(q_e)^{-1} = \operatorname{Ad}_h \circ
\widehat{\mathbb{I}}(q_e)^{-1} \circ \operatorname{Ad}^\ast_h$ for
any $h \in G_{q_e}$. Thus, for $h \in G_{q_e}$, the second summand
in $A_{ab}$ becomes
\begin{align*}
&\left\langle \left(T_{q_e}\mathbb {I}(h\cdot v_{q_e}) \circ
\widehat{\mathbb {I}}(q_e)^{-1} \circ T_{q_e} \overset{\thicksim}
{\mathbb {I}}(h\cdot v_{q_e})\right) \xi_a, \xi_b
\right\rangle \\
&= \left\langle \left(\operatorname{Ad}^\ast_{h^{-1}} \circ
T_{q_e} \mathbb{I}( v_{q_e}) \circ \operatorname{Ad}_{h^{-1}}
\circ \widehat{\mathbb {I}}(q_e)^{-1} \circ
\operatorname{Ad}^\ast_{h^{-1}}\circ
T_{q_e}\overset{\thicksim}{\mathbb {I}}(v_{q_e}) \right)
\xi_a,
\xi_b \right\rangle\\
&= \left\langle \left(\operatorname{Ad}^\ast_{h^{-1}} \circ
T_{q_e} \mathbb{I}( v_{q_e})
\circ \widehat{\mathbb {I}}(q_e)^{-1} \circ
T_{q_e}\overset{\thicksim}{\mathbb {I}}(v_{q_e}) \right)
\xi_a,
\xi_b \right\rangle\\
&= \left\langle \left(
T_{q_e} \mathbb{I}( v_{q_e})
\circ \widehat{\mathbb {I}}(q_e)^{-1} \circ
T_{q_e}\overset{\thicksim}{\mathbb {I}}(v_{q_e}) \right)
\xi_a,
\operatorname{Ad}_{h^{-1}}\xi_b \right\rangle\\
&= \left\langle \left(
T_{q_e} \mathbb{I}( v_{q_e})
\circ \widehat{\mathbb {I}}(q_e)^{-1} \circ
T_{q_e}\overset{\thicksim}{\mathbb {I}}(v_{q_e}) \right)
\xi_a,
\xi_b \right\rangle
\end{align*}
since $\operatorname{Ad}_{h^{-1}}\xi_b = 0 $ because $h \in
G_{q_e}$ and $\xi_b \in \mathfrak{g}_{q_e}$. This shows that the
second summand in $A_{ab}$ is $G_{q_e}$-- invariant.
Next, we show that the first summand in $A_{ab}$ is $G_{q_e}$--
invariant. To see this note that
\[
\left\langle A_{v_{q_e}}\xi_a, \xi_b \right\rangle
= \left.\frac{\partial}{\partial
\tau}\right|_{\tau= 0}\left\langle T_{\tau v_{q_e}}(\mathbb {I}
\circ \operatorname{Exp}_{q_e})(v_{q_e})\xi_a, \xi_b \right\rangle
= \left.\frac{\partial^2}{\partial
\tau^2}\right|_{\tau= 0} \left\langle
\mathbb{I}(\operatorname{Exp}_{q_e}(\tau v_{q_e}))\xi_a, \xi_b
\right\rangle.
\]
Therefore, for any $h \in G_{q_e} $ we get from
\eqref{equivariance of I}
\begin{align*}
\left\langle A_{h\cdot v_{q_e}}\xi_a, \xi_b \right\rangle
&= \left.\frac{\partial^2}{\partial
\tau^2}\right|_{\tau= 0} \left\langle
\mathbb{I}(\operatorname{Exp}_{q_e}(\tau h \cdot v_{q_e}))\xi_a,
\xi_b
\right\rangle
= \left.\frac{\partial^2}{\partial
\tau^2}\right|_{\tau= 0} \left\langle
\mathbb{I}(h \cdot \operatorname{Exp}_{q_e}(\tau  v_{q_e}))\xi_a,
\xi_b \right\rangle \\
&= \left.\frac{\partial^2}{\partial
\tau^2}\right|_{\tau= 0} \left\langle
\operatorname{Ad}_{h^{-1}}^\ast
\mathbb{I}(\operatorname{Exp}_{q_e}(\tau  v_{q_e}))
\operatorname{Ad}_{h^{-1}}\xi_a,
\xi_b \right\rangle\\
&= \left.\frac{\partial^2}{\partial
\tau^2}\right|_{\tau= 0} \left\langle
\mathbb{I}(\operatorname{Exp}_{q_e}(\tau  v_{q_e}))
\operatorname{Ad}_{h^{-1}}\xi_a,
\operatorname{Ad}_{h^{-1}}\xi_b \right\rangle\\
&= \left.\frac{\partial^2}{\partial
\tau^2}\right|_{\tau= 0} \left\langle
\mathbb{I}(\operatorname{Exp}_{q_e}(\tau  v_{q_e}))\xi_a,
\xi_b \right\rangle
= \left\langle A_{v_{q_e}}\xi_a, \xi_b \right\rangle,
\end{align*}
as required.
\end{proof}
\begin{proposition}
\label{determination of xi} The equation $\varphi (\tau
,v_{q_{e}},\mu ,\xi)=0 $ for $(\tau ,v_{q_{e}},\mu ,\xi) \in
I\times \left(B\cap (T_{q_{e}}Q)_{\{e\}}\setminus
\mathcal{Z}_{\mu}\right) \times \mathfrak{g}^{\ast }\times
\mathfrak{g}_{q_{e}}$ has a unique smooth solution $\xi (\tau
,v_{q_{e}},\mu )\in \mathfrak{g}_{q_{e}}$ for $(\tau, v_{q_e},
\mu) \in I\times (B\cap (T_{q_{e}}Q)_{\{e\}}\setminus
\mathcal{Z}_{\mu})\times \mathfrak{g}^{\ast }$.
\end{proposition}
\begin{proof} Denote by $D_\xi$ the Fr\'echet derivative relative
to the variable $\xi \in \mathfrak{g}_{q_e}$. Recall that
$\xi_0(v_{q_e}, \mu) \in \mathfrak{g}_{q_e}$ is the unique
solution of the equation $\frac{\partial ^{2} \varphi}{\partial
\tau ^{2}}(0,v_{q_{e}},\mu ,\xi )=0$. Formulas
\eqref{second tau derivative of phi at zero} and \eqref{first
tau derivative of eta at zero} yield
\begin{align}
\label{expanded second tau derivative of phi at zero}
\frac{\partial ^{2}\varphi }{\partial \tau^{2}} (0,v_{q_{e}},\mu,
\xi )
&= (Id-\Pi)\left[ A_{v_{q_e}}(\xi +\eta _{\mu})
- 2\left(T_{q_e}\mathbb{I}(v_{q_{e}}) \circ
\widehat{\mathbb{I}}(q_e)^{-1} \circ T_{q_e}\overset{\thicksim}
{\mathbb {I}}(v_{q_e})\right) \xi \right. \\
& \phantom{T_{q_e}\overset{\thicksim}
{\mathbb{I}}(v_{q_e})}
- 2\left(T_{q_e}\mathbb{I}(v_{q_{e}}) \circ
\widehat{\mathbb {I}}(q_e)^{-1} \circ
T_{q_e}\widehat{\mathbb{I}}(v_{q_e}) \circ
\widehat{\mathbb{I}} (q_{e})^{-1}\right)(\Pi _{1}\mu )
\nonumber \\
&\left. \phantom{T_{q_e}\overset{\thicksim}
{\mathbb{I}}(v_{q_e})}
+ 2\left(T_{q_e}\mathbb{I}(v_{q_{e}}) \circ\widehat{\mathbb{I}}
(q_{e})^{-1}\right)(\beta'(\mu))  - 2\beta ''(\mu) \right]
\nonumber
\end{align}
and hence
\[
D_\xi\frac{\partial ^{2}\varphi }{\partial \tau^{2}} (0,v_{q_{e}},
\mu, \xi_0(v_{q_e}, \mu))
= (Id-\Pi)\left[ A_{v_{q_e}}|_{\mathfrak{g}_{q_e}}
- 2T_{q_e}\mathbb{I}(v_{q_{e}}) \circ
\widehat{\mathbb{I}}(q_e)^{-1} \circ T_{q_e}\overset{\thicksim}
{\mathbb {I}}(v_{q_e})
\right] : \mathfrak{g}_{q_e} \rightarrow \mathfrak{k}^\circ.
\]
We shall prove that this linear map is injective. To see this,
note that relative to the basis $\{\xi_1, \dots , \xi_p\}$ of
$\mathfrak{g}_{q_e}$ this linear operator has matrix $A $ by
\eqref{A}. Thus, if $v_{q_e} \notin \mathcal{Z}_{\mu}$, this
matrix is invertible. In particular, this linear operator is
injective. 
 
Since $\mathfrak{g} = \mathfrak{g}_{q_e} \oplus
\mathfrak{k}$, it follows that $\dim \mathfrak{g}_{q_e} = \dim
\mathfrak{g} - \dim \mathfrak{k} = \dim \mathfrak{k}^\circ $, so
the injectivity of the map $D_\xi\frac{\partial ^{2}\varphi
}{\partial \tau^{2}} (0,v_{q_{e}}^0, \mu^0, \xi_0(v_{q_e}^0,
\mu^0))$ implies that it is an isomorphism. Therefore, if
$v_{q_{e}}\in B\cap (T_{q_{e}}Q)_{\{e\}}\setminus
\mathcal{Z}_{\mu}$ is near $v_{q_e}^0 $, the implicit function
theorem, guarantees the existence of an open neighborhood
$V_{0}\subset I\times (B\cap (T_{q_{e}}Q)_{\{e\}}\setminus
\mathcal{Z}_{\mu})\times
\mathfrak{g}^{\ast }$ containing $(0, v_{q_e}^0, \mu^0) \in
\{0\}\times (B\cap (T_{q_{e}}Q)_{\{e\}}\setminus
\mathcal{Z}_{\mu}) \times
\mathfrak{g}^{\ast }$ and of a unique smooth function $\xi
:V_{0}\rightarrow \mathfrak{g}_{q_{e}}$ satisfying $\varphi (\tau
,v_{q_{e}},\mu , \xi(\tau
,v_{q_{e}},\mu) )=0$ such that $\xi(0,v_{q_{e}}^0,\mu^0) =
\xi_0(v_{q_e}^0,
\mu^0)$.
On the other hand, for $\tau \neq 0$, the equation $\varphi (\tau
,v_{q_{e}},\mu ,\cdot )=0$ has a unique solution for $\xi $,
namely the $\mathfrak{g}_{q_{e}}$-component of
$\mathbb{I}(\operatorname{Exp} _{q_{q}}(\tau
v_{q_{e}}))^{-1}\beta (\tau ,\mu )$, which is a smooth function
of $\tau ,v_{q_{e}},\mu $. This is true since $\xi + \eta =
\mathbb{I}(\operatorname{Exp} _{q_{q}}(\tau v_{q_{e}}))^{-1}\beta
(\tau ,\mu )$ by construction and we determined the two
components $\xi\in \mathfrak{g}_{q_e} $ and $\eta \in
\mathfrak{k}$ in $\mathfrak{g} = \mathfrak{g}_{q_e} \oplus
\mathfrak{k}$ via the Lyapunov-Schmidt method, precisely in order
that this equality be satisfied. Therefore, the solution
$\xi(\tau, v_{q_e}, \mu)$ obtained above by the implicit
function theorem must coincide with the $\mathfrak{g}_{q_{e}}$-component of
$\mathbb{I}(\operatorname{Exp} _{q_{q}}(\tau
v_{q_{e}}))^{-1}\beta (\tau ,\mu )$ for $\tau>0 $.
Since this entire argument involving the Lyapunov-Schmidt
procedure was carried out for any $(v_{q_e}^0, \mu^0)$, it
follows that the equation
$\varphi (\tau,v_{q_{e}},\mu ,\xi)=0$ has a unique smooth
solution $\xi (\tau ,v_{q_{e}},\mu )\in \mathfrak{g}_{q_{e}}$
for $(\tau, v_{q_e},
\mu) \in I\times (B\cap (T_{q_{e}}Q)_{\{e\}}\setminus
\mathcal{Z}_{\mu})\times \mathfrak{g}^{\ast }$.
\end{proof}
\begin{remark}
\label{remark on smoothness}
\normalfont
The previous proposition says that if we define
\begin{equation*}
\zeta(\tau ,v_{q_{e}},\mu )=\mathbb{I}(\operatorname{Exp}
_{q_{e}}(\tau v_{q_{e}}))^{-1}\beta (\tau ,\mu )
\end{equation*}
on $(I\setminus \{0\})\times (B\cap (T_{q_{e}}Q)_{\{e\}}\setminus
\mathcal{Z}_{\mu})\times \mathfrak{g}^{\ast }$, then $\zeta(\tau
,v_{q_{e}},\mu )$ can be smoothly extended for $\tau=0$.
We have, in fact, $\zeta(\tau ,v_{q_{e}},\mu )= \xi(\tau
,v_{q_{e}},\mu ) + \eta (\tau ,v_{q_{e}},\mu ,\xi(\tau
,v_{q_{e}},\mu ) )$, where $\eta(\tau ,v_{q_{e}},\mu, \xi )$ was
found in the first step of the Lyapunov-Schmidt procedure and
$\xi(\tau ,v_{q_{e}},\mu )$ in the second step, as given in
Proposition \ref{determination of xi}. Note also that $\zeta(0,
v_{q_e}, \mu) = \xi_0(v_{q_e}, \mu) + \widehat{\mathbb
{I}}(q_e)^{-1} \Pi_1 \mu \in \mathfrak{t}$.
\end{remark}
\subsection{A simplified version of the amended potential
criterion} At this point we have  a candidate for a bifurcating
branch from the set of relative equilibria $\mathfrak{t} \cdot
q_e $. This branch will start at $\zeta(0, v_{q_e}, \mu)_Q(q_e)
\in \mathfrak{t}\cdot q_e \subset T_{q_e}Q$. By Lemma \ref{same
isotropy}, the isotropy subgroup of $\zeta(0, v_{q_e}, \mu)_Q(q_e)
$ equals $G_{q_e} $, for any $v_{q_e} \in B\cap
(T_{q_{e}}Q)_{\{e\}}\setminus \mathcal{Z}_{\mu}$ and $\mu\in
\mathfrak{g}^\ast$. The isotropy groups of the points on the
curve $\zeta(\tau, v_{q_e}, \mu)_Q(\operatorname{Exp}_{q_{e}}
(\tau v_{q_e}))$, for $\tau\neq 0 $, are all trivial, by
construction. Hence $\zeta(\tau, v_{q_e},
\mu)_Q(\operatorname{Exp}_{q_{e}}(\tau v_{q_e}))$ is a curve
that has the properties of the bifurcating branch of relative
equilibria with broken symmetry that we are looking for. We do
not know yet that all points on this curve are in fact relative
equilibria. Thus, we shall search for conditions on $v_{q_e}$ and
$\mu$ that guarantee that each point on  the curve
$\tau\mapsto\zeta(\tau, v_{q_e},
\mu)_Q(\operatorname{Exp}_{q_{e}}(\tau v_{q_e}))$ is a relative
equilibrium. This will be done by using the amended potential
criterion (see Proposition \ref{amended potential criterion}) which
is applicable because all base points of this curve, namely
$\operatorname{Exp}_{q_e}(\tau v_{q_e})$, have trivial isotropy
for $\tau\neq 0$. To carry this out, we need some additional
geometric information.
\medskip
From standard theory of proper Lie group actions (see e.g.
\cite{dk}, \S 2.3, or \cite{kawakubo}) it follows that
the map
\begin{equation}
\label{corespondenta} 
[v_{q_{e}},\mu ]_{G_{q_{e}}} \in (B
\times \mathfrak{g}^{\ast })/G_{q_{e}}
\longmapsto [\operatorname{Exp}_{q_{e}}(v_{q_{e}}),\mu ]_{G} \in
((G \cdot \operatorname{Exp}_{q_e} B) \times \mathfrak{g}^{\ast
})/G
\end{equation}
is a homeomorphism of $(B\times \mathfrak{g}^{\ast })/G_{q_{e}}$
with $((G \cdot \operatorname{Exp}_{q_e} B) \times
\mathfrak{g}^{\ast })/G$ and that its restriction to $ ((B\cap
(T_{q_{e}}Q)_{\{e\}} \setminus \mathcal{Z}_{\mu})\times
\mathfrak{g}^{\ast })/G_{q_{e}}$ is a diffeomorphism onto its
image. We think of a pair $(\operatorname{Exp}_{q_e}(v_{q_e}),
\mu)$ as the base point of a relative equilibrium and its
momentum value. All these relative equilibria come in
$G$-orbits. The homeomorphism \eqref{corespondenta} allows the
identification of $G$-orbits of relative equilibria with
$G_{q_e}$-orbits of certain pairs $(v_{q_e},\mu)$. We shall work
in what follows on both sides of this identification, based on
convenience.
 We will need the
following lemma, which is a special case of stability of
the transversality of smooth maps (see e.g. \cite{gp}).
\begin{lemma}
Let $G$ be a Lie group acting on a Riemannian manifold $Q$, $q\in
Q$, and let $\mathfrak{k}\subset \mathfrak{g}$ be a subspace
satisfying $\mathfrak{k}\cap \mathfrak{g}_{q}=\{0\}$. Let
$V\subset T_{q}Q$ be a subspace such that $\mathfrak{k}\cdot
q\oplus V=T_{q}Q$. Then there is an $\epsilon >0$ such that if
$\|v_{q}\| <\epsilon $,
\begin{equation*}
T_{\operatorname{Exp} _{q}(v_{q})}Q=\mathfrak{k}\cdot \operatorname{Exp} _{q}(v_{q})\oplus
(T_{v_{q}}\operatorname{Exp} _{q})V.
\end{equation*}
\end{lemma}
To deal with $G $-orbits  of relative equilibria, we need a
different splitting of the same nature. The following result
is modeled on a proposition in \cite{hm}.
\begin{proposition}
\label{decomposition of the tangent space adapted to sigma bar}
Let $v_{q_e} \in B\cap (T_{q_{e}}Q)_{\{e\}} \setminus
\mathcal{Z}_{\mu}$ be given. Consider the principal
$G_{q_e}$-bundle $B\cap (T_{q_{e}}Q)_{\{e\}} \setminus
\mathcal{Z}_{\mu} \rightarrow [B\cap (T_{q_{e}}Q)_{\{e\}}
\setminus \mathcal{Z}_{\mu}]/G_{q_{e}}$ (this is implied by Lemma
\ref {properties of Z}). Let $\widetilde{U}$ be a neighborhood of
$[0_{q_e}] \in (T_{q_e}Q)/G_{q_e}$ and  define the open set $U: =
\widetilde{U} \cap [B\cap (T_{q_{e}}Q)_{\{e\}} \setminus
\mathcal{Z}_{\mu}]/G_{q_{e}} $ in  $[B\cap (T_{q_{e}}Q)_{\{e\}}
\setminus \mathcal{Z}_{\mu}]/G_{q_{e}} $. Let $\sigma :U\subset
[B\cap (T_{q_{e}}Q)_{\{e\}})\setminus
\mathcal{Z}_{\mu}]/G_{q_{e}}\rightarrow  B\cap
(T_{q_{e}}Q)_{\{e\}}\setminus \mathcal{Z}_{\mu}$ be a smooth
section, $[v_{_{q_{e}}}]\in U$, and $\overline{\sigma
}:=\operatorname{Exp} _{q_{e}}\circ \sigma : U \rightarrow Q$.
Then there exists $\epsilon >0$ such that for $ 0<\tau <\epsilon
$ sufficiently small, we have
\begin{equation*}
T_{\overline{\sigma }([\tau v_{q_{e}}])}Q
=\mathfrak{t}\cdot \overline{\sigma }([\tau v_{q_{e}}])
\oplus T_{[\tau v_{q_{e}}]}
\overline{\sigma}(T_{[\tau v_{q_{e}}]}U)
\oplus (T_{\sigma ([\tau v_{q_{e}}])}\operatorname{Exp}_{q_{e}})
(\mathfrak{k} _{2}\cdot q_{e}).
\end{equation*}
\end{proposition}
\begin{proof}
Since $\mathfrak{g}= \mathfrak{k}_0 \oplus \mathfrak{k}_1 \oplus
\mathfrak{k}_2 $ and $\mathfrak{k}_0 = \mathfrak{g}_{q_e}$ we
have
$T_{q_{e}}Q=\mathfrak{k}_{1}\cdot q_{e}\oplus
\mathfrak{k}_{2}\cdot q_{e}\oplus (\mathfrak{g}\cdot
q_{e})^{\perp }$. Apply the above lemma with
$\mathfrak{k}=\mathfrak{k}_{1}$ and $V=\mathfrak{k}_{2}\cdot
q_{e}\oplus (\mathfrak{g}\cdot q_{e})^{\perp }$. For the
$\epsilon >0 $ in the statement choose $\tau $ such that  $0<
\tau< \epsilon $ and $\| \sigma([ \tau v_{q_e}]) \| <
\epsilon$. Then
\begin{align}
\label{decomposition of q for the section}
T_{\overline{\sigma }([\tau v_{q_{e}}])}Q
&=\mathfrak{k}_{1}\cdot \overline{
\sigma }([\tau v_{q_{e}}])
\oplus (T_{\sigma ([\tau v_{q_{e}}])}\operatorname{Exp}
_{q_{e}})(\mathfrak{k}_{2}\cdot q_{e} \oplus (\mathfrak{g}\cdot
q_e) ^\perp)\\
&= \mathfrak{k}_{1}\cdot \overline{\sigma }([\tau v_{q_{e}}])
\oplus (T_{\sigma ([\tau v_{q_{e}}])}\operatorname{Exp}
_{q_{e}})((\mathfrak{g}\cdot q_e) ^\perp )
\oplus (T_{\sigma ([\tau v_{q_{e}}])}\operatorname{Exp}
_{q_{e}})(\mathfrak{k}_2 \cdot q_e). \nonumber
\end{align}
since $\operatorname{Exp}_{q_e}$ is a diffeomorphism on $B
\subset (\mathfrak{g}\cdot q_e)^\perp$. Since $(\sigma ,U)$ is a
smooth local section, $\mathcal{Z}_{\mu}$ is closed and
$G_{q_e}$-invariant in $B\cap (T_{q_e}Q)_{\{e\}}$, and
$(T_{q_e}Q)_{\{e\}}$ is open in $T_{q_e} Q $, it follows that
$B\cap (T_{q_e}Q)_{\{e\}}$ is open in $(\mathfrak{g}\cdot
q_e)^\perp$ and thus  we get
\begin{equation*}
(\mathfrak{g}\cdot q_e)^\perp = T_{\sigma ([\tau
v_{q_{e}}])}(B\cap (T_{q_{e}}Q)_{\{e\}}\setminus
\mathcal{Z}_{\mu}) = T_{[\tau v_{q_{e}}]}\sigma (T_{[\tau
v_{q_{e}}]}U) \oplus \mathfrak{k}_{0}\cdot \sigma ([\tau
v_{q_{e}}]),
\end{equation*}
where $\mathfrak{k}_{0}\cdot \sigma ([\tau v_{q_{e}}]) =
\{\zeta_{T_{q_e}Q} (\sigma([ \tau v_{q_e}]) \mid \zeta\in
\mathfrak{k}_0 \}$.  The
$G_{q_e}$-equivariance of
$\operatorname{Exp}_{q_e}$ implies that
\[
T_{u_{q_{e}}}\operatorname{Exp} _{q_{e}}(\xi
_{T_{q_e}Q}(u_{q_{e}}))=\xi _{Q}(\operatorname{Exp}
_{q_{e}}(u_{q_{e}})) \quad \text{for~all} \quad
\xi \in \mathfrak{k}_{0},  \quad u_{q_{e}}\in T_{q_{e}}Q
\]
 and hence
\begin{align}
\label{decomposition of one summand}
&(T_{\sigma ([\tau v_{q_{e}}])}\operatorname{Exp}
_{q_{e}})((\mathfrak{g}\cdot q_e) ^\perp ) \\
&\qquad \qquad =(T_{\sigma([\tau v_{q_{e}}])}\operatorname{Exp}
_{q_{e}}
\circ T_{[\tau v_{q_e}]}\sigma )(T_{[\tau v_{q_{e}}]}U)\oplus
(T_{\sigma ([\tau v_{q_{e}}])}\operatorname{Exp}
_{q_{e}})(\mathfrak{k}_{0}\cdot
\sigma ([\tau v_{q_{e}}])) \nonumber  \\
&\qquad \qquad  =T_{[\tau v_{q_{e}}]}\overline{\sigma
}(T_{[\tau v_{q_{e}}]}U)\oplus
\mathfrak{k}_{0}\cdot \overline{\sigma }([\tau v_{q_{e}}]).
\nonumber
\end{align}
Introducing \eqref{decomposition of one summand} in
\eqref{decomposition of q for the section} and taking into
account that $\mathfrak{t}=\mathfrak{k}_{0}\oplus
\mathfrak{k}_{1}$ we get the statement of the proposition.
\end{proof}
We want to find pairs $(v_{q_e}, \mu) $ such that
$\mathbf{d}V_{\beta(\tau, \mu)}  (\operatorname{Exp} _{q_{e}}(\tau
v_{q_{e}}))=0$ for $\tau> 0 $. Since $V_{\beta(\tau, \mu)}$ is
$G_{\beta(\tau,
\mu)}$-invariant, this condition will hold if we only verify it
on a subspace of $T_{\operatorname{Exp}_{q_e}(\tau v_{q_e})} Q$
complementary to $\mathfrak{g}_{\beta(\tau,\mu)} \cdot
\operatorname{Exp} _{q_{e}}(\tau v_{q_{e}}) = \mathfrak{t}\cdot
\operatorname{Exp} _{q_{e}}(\tau v_{q_{e}})$. The previous
decomposition of the tangent space immediately yields the
following result.
\begin{corollary}
\label{criteriu}
Suppose that $\mu\in \mathfrak{g}^\ast$ is such
that $\mathfrak{g}_{\beta( \tau, \mu)} = \mathfrak{t}$ for all
$\tau$ in a neighborhood of zero.  Let $U$ and $\sigma $ be as in
Proposition  \ref{decomposition of the tangent space adapted to
sigma bar}, $[v_{q_{e}}]\in U$, and $\overline{\sigma
}:=\operatorname{Exp} _{q_{e}}\circ \sigma $. Then there is an
$\epsilon >0$ such that $\mathbf{d}V_{\beta (\tau,
\mu)}(\overline{\sigma }([\tau v_{q_{e}}])=0$ if and only if
$\mathbf{d}(V_{\beta (\tau, \mu )}\circ \overline{\sigma })([\tau
v_{q_{e}}])=0$ and $\mathbf{d}(V_{\beta (\tau, \mu )}\circ
\operatorname{Exp} _{q_{e}})( \sigma ([\tau v_{q_{e}}]))|
_{\mathfrak{k}_{2}\cdot q_{e}}=0$ for $ 0<\tau <\epsilon $.
\end{corollary}
\subsection{The study of two auxiliary functions}
Let $I $ be an open interval containing zero. Recall that $p =
\dim \mathfrak{g}_{q_e} = \dim \mathfrak{m}_0 $. Let
$\vartheta_{1}$ be an element of a basis
$\{\vartheta _{1},\vartheta _{2},...,
\vartheta _{p}\}$ for $\mathfrak{m}_{0}$ and define $\beta
:(I\setminus \{0\})\times (\mathfrak{m}_{1}\oplus
\mathfrak{m}_{2})\rightarrow
\mathfrak{g}^{\ast }$ by
\begin{equation*}
\beta (\tau ,\mu)=\Pi _{1}\mu+\tau \Pi_{2}\mu
+\tau ^{2}\vartheta _{1},
\end{equation*}
where $\Pi_1 : \mathfrak{g}^\ast\rightarrow \mathfrak{m}_1 =
\mathbb {I}(q_e) \mathfrak{t}$ and $\Pi_2: \mathfrak{g}^\ast
\rightarrow \mathfrak{m}_2 = \mathfrak{t}^\circ $. Notice
that this function is a particular case of
\begin{equation*}
\beta (\tau ,\mu )=\Pi _{1}\mu +\tau \beta ^{\prime }(\mu )+\tau ^{2}\beta
^{\prime \prime }(\mu ),
\end{equation*}
by choosing $\beta' (\mu) = \Pi_2 \mu$ and $\beta '' (\mu) =
\vartheta_1 $.
Recall that $\mathbb {I}(q_e) = \mathfrak{m}_1
\oplus\mathfrak{m}_2 $ by Lemma \ref{moment of  inertia
isomorphism} and that $\mathbf{J}_L(\mathfrak{g}\cdot q_e) =
\mathbb {I}(q_e) \mathfrak{g}$ from the definition of
$\mathbf{J}_L $.
\begin{theorem}
\label{theorem about F}
The smooth function $F_{1}:(I\setminus \{0\})\times U\times
\mathbf{J}_{L}(\mathfrak{g}\cdot q_{e})\rightarrow \mathbb{R}$
defined by
\begin{equation*}
F_{1}(\tau ,[v_{q_{e}}],\mu ):=(V_{\beta (\tau, \mu )}\circ
\overline{\sigma })(\tau [v_{q_{e}}]).
\end{equation*}
can be extended to a smooth function on $I\times U\times
\mathbf{J}_{L}(\mathfrak{g}\cdot q_{e})$, also denoted by $F_1
$. In addition
\begin{equation*}
F_{1}(\tau ,[v_{q_{e}}],\mu )=F_{0}(\mu )+\tau ^{2}F(\tau
,[v_{q_{e}}],\mu ).
\end{equation*}
where $F_{0}$, $F$ are defined on $\mathbf{J}_{L}(
\mathfrak{g}\cdot q_{e})$ and on $I\times U\times \mathbf{J}_{L}(
\mathfrak{g}\cdot q_{e})$ respectively.
\end{theorem}
\begin{proof}
Denote $v_{q_{e}}:=\sigma ([v_{q_{e}}])\in B\cap
(T_{q_{e}}Q)_{\{e\}} \setminus \mathcal{Z}_{\mu}$. One can easily
see that
\begin{equation*}
(V_{\beta (\tau, \mu)}\circ \overline{\sigma })(\tau \lbrack
v_{q_{e}}])=V(\operatorname{Exp} _{q_{e}}(\tau
v_{q_{e}}))+\frac{1}{2}\left\langle\beta (\tau ,\mu
),\mathbb{I}(\operatorname{Exp}_{q_{e}}(\tau v_{q_{e}}))^{-1}
\beta(\tau ,\mu )\right\rangle.
\end{equation*}
By Remark \ref{remark on smoothness}, the second term is smooth
even in a a neighborhood of $\tau= 0 $. Since the first term is
obviously smooth, it follows that $V_{\beta (\tau, \mu)}\circ
\overline{\sigma }$ is smooth also in a neighborhood of $\tau= 0
$. This is the smooth extension of $F_1 $ in the statement.
Let $\{\xi _{1},...,\xi _{p}\}$ be a
basis for $\mathfrak{g}_{q_{e}} \subset \mathfrak{t}$. Then,
again by Remark \ref{remark on smoothness}, we have
\begin{equation*}
\mathbb{I}(\operatorname{Exp}_{q_{e}}(\tau v_{q_{e}}))^{-1}\beta
(\tau ,\mu )=\overset{p}{\underset{a=1}{
\sum }}\alpha _{a}(\tau ,v_{q_{e}},\mu )\xi _{a}+\eta \left(\tau
,v_{q_{e}},\mu ,
\overset{p}{\underset{a=1}{\sum }}\alpha _{a}(\tau ,v_{q_{e}},\mu
)\xi _{a}\right)
\end{equation*}
where $\alpha _{1},...,\alpha _{p},\eta $ are smooth real
functions of all their arguments. In what follows we will denote
\begin{equation*}
\eta \left(\tau ,v_{q_{e}},\mu ,\overset{p}{\underset{a=1}{\sum
}}\alpha _{a}(\tau ,v_{q_{e}},\mu )\xi _{a}\right)
=\eta (\tau,v_{q_{e}},\mu ,\alpha _{1}(\tau ,v_{q_{e}},\mu
),...,\alpha _{p}(\tau ,v_{q_{e}},\mu )).
\end{equation*}
Let $\mu \in \mathbf{J}_{L}(\mathfrak{g}\cdot q_{e})=
\mathfrak{m}_1 \oplus\mathfrak{m}_2$ and $v_{q_{e}}\in B\cap
(T_{q_{e}}Q)_{\{e\}}\backslash \mathcal{Z}_{\mu}$. Since  in the
computations that follow, the arguments $v_{q_e}$ and $\mu$ play
the role of parameters, we shall denote temporarily  $\alpha
_{a}(\tau )=\alpha _{a}(\tau, v_{q_{e}},\mu)$, $a\in
\{1,...,p\}$, and $\eta (\tau ,\alpha _{1},...,\alpha _{p})=\eta
(\tau ,v_{q_{e}},\mu ,\alpha _{1}(\tau ,v_{q_{e}},\mu ),...,\alpha
_{p}(\tau ,v_{q_{e}},\mu ))$. Then by \eqref{first tau derivative
of eta at zero} we get
\begin{align*}
\frac{\partial \eta }{\partial \tau }(0,\alpha _{1},...,
\alpha _{p})
= &-\overset{p}{\underset{a=1}{\sum }}\alpha_{a}
\left(\widehat{\mathbb{I}}(q_{e})^{-1}\circ
T_{q_e}\overset{\thicksim}{\mathbb{I}}
(v_{q_{e}})\right)\xi _{a} \\
& -\left(\widehat{\mathbb{I}}(q_{e})^{-1} \circ
T_{q_e}\widehat{\mathbb{I}}(v_{q_{e}}) \circ
\widehat{\mathbb{I}}(q_{e})^{-1} \right)\Pi_{1}\mu
+\widehat{\mathbb{I}}(q_{e})^{-1}\Pi_{2}\mu.
\end{align*}
Formula \eqref{pi composed with Phi} shows that
\[
\frac{\partial \eta }{\partial \alpha_a }(0,\alpha_1, \dots ,
\alpha_p)  = 0
\]
Note that
\begin{equation*}
\left.V_{\beta (\tau ,\mu )}(\operatorname{Exp}_{q_{e}}(\tau
v_{q_{e}}))\right|_{\tau =0} = V(q_{e})+\frac{1}{2}
\left\langle \Pi_{1}\mu,\widehat
{\mathbb{I}}(q_{e})^{-1}\Pi_{1}\mu\right\rangle
\end{equation*}
is independent of $v_{q_{e}}$. This shows that $F_1(0, [v_{q_e}],
\mu) = F_0(\mu)$ for some smooth function on $\mathfrak{m}_1
\oplus \mathfrak{m}_2 $.
Using Remark \ref{remark on smoothness}, we get
\begin{align*}
\left.\frac{d}{d\tau }\right|_{\tau =0}&V_{\beta
(\tau ,\mu )}(\operatorname{Exp}_{q_{e}}(\tau v_{q_{e}}))
=\mathbf{d}V(q_{e})(v_{q_{e}}) + \frac{1}{2}\left\langle
\Pi_{2}\mu ,\overset{p}{
\underset{a=1}{\sum }}\alpha _{a}(0)\xi _{a}+\eta (0,\alpha _{1},...,\alpha
_{p})\right\rangle\\
&+\frac{1}{2}\left\langle \Pi_1\mu, \sum_{a=1}^p \frac{\partial
\alpha_a}{\partial \tau}(0)\left(\xi_a + \frac{\partial
\eta}{\partial \alpha_a}(0, \alpha_1, \dots , \alpha_p)   \right)
+ \frac{\partial \eta}{\partial \tau}(0, \alpha _1, \dots ,
\alpha_p)  \right\rangle.
\end{align*}
The first term
$\mathbf{d}V(q_{e})=0$ by 
Proposition \ref{montaldi} (i).
Since $\eta(0, v_{q_e}, \mu, \xi) =
\eta_\mu = \widehat{\mathbb {I}}(q_e)^{-1}\Pi_1 \mu \in
\mathfrak{t}$ by Proposition \ref{belongs}, we get
\[
\overset{p}{\underset{ a=1}{\sum
}}\alpha _{a}(0)\xi _{a}+\eta (0,\alpha _{1},...,\alpha _{p})=
\overset{p}{\underset{a=1}{\sum }}\alpha _{a}(0)\xi _{a}
+\widehat{\mathbb{I}}(q_{e})^{-1}\Pi_{1}\mu\in \mathfrak{t}.
\]
Thus the second term vanishes because
$\mathfrak{m}_{2}=\mathfrak{t}^{\circ }$. 
As $\frac{\partial \eta }{\partial \alpha _{a}}(0, \alpha_1,
\dots , \alpha_p) =0$ and
$\mathfrak{m}_{1}$ annihilates $\mathfrak{g} _{q_{e}}$, the third
term becomes
\begin{align*}
\left\langle \Pi_{1}\mu,\frac{\partial \eta }{\partial \tau }
(0,\alpha _{1},...,\alpha_{p})\right\rangle
=&-\overset{p}{\underset{a=1}{\sum}}\alpha _{a}
\left\langle \Pi_1\mu,
\left(\widehat{\mathbb{I}}(q_{e})^{-1}\circ
T_{q_e}\overset{\thicksim}{\mathbb{I}}
(v_{q_{e}})\right)\xi _{a}\right\rangle \\
&-\left\langle\Pi_{1}\mu,\left(\widehat{\mathbb{I}}(q_{e})^{-1}
\circ T_{q_e}\widehat{\mathbb{I}}(v_{q_{e}}) \circ
\widehat{\mathbb{I}}(q_{e})^{-1} \right)\Pi_{1}\mu\right\rangle\\
&+\left\langle\Pi_{1}\mu,\widehat{\mathbb{I}}(q_{e})^{-1}\Pi_{2}\mu
\right\rangle.
\end{align*}
We will prove that each summand in this expression vanishes.
$\bullet$ Since
$\langle \mathfrak{m}_0 , \mathfrak{k}_1 \rangle = 0$, we get
\begin{align*}
&\left\langle \Pi_1\mu,
\left(\widehat{\mathbb{I}}(q_{e})^{-1}\circ
T_{q_e}\overset{\thicksim}{\mathbb{I}}
(v_{q_{e}})\right)\xi _{a}\right\rangle
= \left\langle
T_{q_e}\overset{\thicksim}{\mathbb{I}}
(v_{q_{e}})\xi_{a},\;
\widehat{\mathbb{I}}(q_{e})^{-1}\Pi_1\mu \right\rangle \\
&\qquad = \left\langle
T_{q_e}{\mathbb{I}}
(v_{q_{e}})\xi_{a},\;
\widehat{\mathbb{I}}(q_{e})^{-1}\Pi_1\mu \right\rangle
= \mathbf{d}\left\langle \mathbb {I}(\cdot) \xi_a,
\eta_\mu \right\rangle(q_e)(v_{q_e}) =0
\end{align*}
by \eqref{differential of i} because $\xi_a \in
\mathfrak{g}_{q_e}$ and $\eta_\mu \in \mathfrak{t}$. Thus the
first summand vanishes.
$\bullet$ The second summand equals
\begin{equation*}
\left\langle\Pi_{1}\mu,\left(\widehat{\mathbb{I}}(q_{e})^{-1}
\circ T_{q_e}\widehat{\mathbb{I}}(v_{q_{e}}) \circ
\widehat{\mathbb{I}}(q_{e})^{-1} \right)\Pi_{1}\mu\right\rangle
= \left\langle T_{q_e}\widehat{\mathbb{I}}(v_{q_{e}})
\eta_\mu, \eta_\mu  \right\rangle
= \left\langle T_{q_e}\mathbb{I}(v_{q_{e}})
\eta_\mu, \eta_\mu  \right\rangle
\end{equation*}
because $\langle \mathfrak{m}_0, \mathfrak{k}_1 \rangle = 0 $. We
shall prove that this term vanishes in the following way. Recall
that $\eta_\mu \in \mathfrak{k}_1 \subset \mathfrak{t}$. For any
$\zeta \in \mathfrak{t}$, hypothesis \textbf{(H)} states that
$\zeta_Q(q_e) $ is a relative equilibrium and thus, by the
augmented potential criterion (see Proposition \ref{augmented
potential criterion}), $\mathbf{d}V_\zeta(q_e) = 0 $. Since
\[
\mathbf{d}V_\zeta(q_e) (u_{q_e}) = \mathbf{d}V(q_e)(u_{q_e}) -
\frac{1}{2} \left\langle T_{q_e} \mathbb{I}(u_{q_e}) \zeta,
\zeta  \right\rangle
\]
for any $u_{q_e} \in T_{q_e} Q $ and $\mathbf{d}V(q_e) = 0 $ by
Proposition \ref{montaldi} (i), it follows that $\left\langle
T_{q_e}
\mathbb{I}(u_{q_e}) \zeta,\zeta\right\rangle = 0 $. Thus the
second summand vanishes.
$\bullet$ The third summand is
\begin{equation*}
\left\langle\Pi_{1}\mu,\widehat{\mathbb{I}}(q_{e})^{-1}\Pi_{2}\mu
\right\rangle
=\langle \Pi_{2}\mu,\eta_\mu \rangle =0
\end{equation*}
because $\mathfrak{m}_{2}=\mathfrak{t}^{\circ }$.
So, we finally conclude that
\begin{equation*}
\left.\frac{d}{d\tau }\right|_{\tau =0}V_{\beta
(\tau ,\mu )}(\operatorname{Exp}_{q_{e}}(\tau v_{q_{e}})) =0
\end{equation*}
and hence, by Taylor's theorem, we have
\begin{equation*}
F_{1}(\tau ,[v_{q_{e}}],\mu )=F_{0}(\mu )+\tau ^{2}F(\tau,
[v_{q_{e}}],\mu )
\end{equation*}
for some smooth function $F $.
\end{proof}
\begin{theorem}
\label{theorem about G}
The smooth function $G_{1}:(I \setminus \{0\}) \times U\times
\mathbf{J}_{L}(\mathfrak{g}\cdot q_{e})\rightarrow
\mathfrak{k}_{2}^{\ast }$  defined by
\begin{equation*}
\left\langle G_{1}(\tau ,[v_{q_{e}}],\mu ),\varsigma
\right\rangle
=\mathbf{d}(V_{\beta (\tau ,\mu )}\circ
\operatorname{Exp}_{q_{e}})(\sigma (\tau \lbrack v_{q_{e}}]))\big(
\varsigma _{Q}(q_{e})\big), \quad \varsigma \in \mathfrak{k}_2,
\end{equation*}
can be smoothly extended to a function on $I\times U\times
\mathbf{J}_{L}(\mathfrak{g}\cdot q_{e})$, also denoted by $G_1$.
In addition,
\begin{equation*}
G_{1}(\tau ,[v_{q_{e}}],\mu )=\tau G(\tau ,[v_{q_{e}}],\mu )
\end{equation*}
where $G:I\times U\times \mathbf{J}_{L}(\mathfrak{g}\cdot
q_{e})\rightarrow \mathfrak{k}_{2}^{\ast }$ is a smooth function.
\end{theorem}
\begin{proof}
We will show that $G_{1}$ is a smooth function at $\tau =0$ and
that $G_{1}(0,[v_{q_{e}}],\mu )=0$.
Let $v_{q_{e}}=\sigma ([v_{q_{e}}])$. Then
\begin{align*}
&\left\langle G_{1}(\tau ,[v_{q_{e}}],\mu ),\varsigma
\right\rangle
= \mathbf{d}V_{\beta(\tau,\mu)}
\left(\operatorname{Exp}_{q_e}(\tau v_{q_e})\right)
\left(T_{\tau v_{q_e}} \operatorname{Exp}_{q_e}
\big(\varsigma_Q(q_e) \big)  \right) \\
&\qquad =\mathbf{d}V(\operatorname{Exp}_{q_{e}}(\tau v_{q_{e}}))
\left(T_{\tau v_{q_e}} \operatorname{Exp}_{q_e}
\big(\varsigma_Q(q_e) \big)  \right) \\
& \qquad \qquad \qquad  +\frac{1}{2}\left\langle \beta (\tau
,\mu),\; T_{\operatorname{Exp}_{q_e}( \tau v_{q_e})}
(\mathbb{I}(\cdot)^{-1})
\left(T_{\tau v_{q_e}} \operatorname{Exp}_{q_e}
\big(\varsigma_Q(q_e) \big)  \right)
\beta (\tau ,\mu )
\right\rangle \\
&\qquad = \mathbf{d}V
\left(\operatorname{Exp}_{q_e}(\tau v_{q_e})\right)
\left(T_{\tau v_{q_e}} \operatorname{Exp}_{q_e}
\big(\varsigma_Q(q_e) \big)  \right)
-\frac{1}{2}\left\langle \beta(\tau,\mu),\;
\phantom{T_{\operatorname{Exp}_{q_e}( \tau v_{q_e})}\mathbb{I}}
\right.\\
&\qquad \qquad \qquad \left.
\left[ \mathbb {I}(\operatorname{Exp}_{q_e}(\tau v_{q_e}))^{-1}
\circ
T_{\operatorname{Exp}_{q_e}( \tau v_{q_e})}\mathbb{I}
\left(T_{\tau v_{q_e}} \operatorname{Exp}_{q_e}
\big(\varsigma_Q(q_e) \big)  \right)
\circ \mathbb {I}(\operatorname{Exp}_{q_e}(\tau v_{q_e}))^{-1}
\right] \beta( \tau, \mu)
\right\rangle \\
&\qquad =\mathbf{d}V
\left(\operatorname{Exp}_{q_e}(\tau v_{q_e})\right)
\left(T_{\tau v_{q_e}} \operatorname{Exp}_{q_e}
\big(\varsigma_Q(q_e) \big)  \right) \\
& \qquad \qquad \qquad  -\frac{1}{2}\left\langle \zeta (\tau,
v_{q_e}, \mu),\; T_{\operatorname{Exp}_{q_e}( \tau v_{q_e})}
\mathbb{I}
\left(T_{\tau v_{q_e}} \operatorname{Exp}_{q_e}
\big(\varsigma_Q(q_e) \big)  \right)
\zeta (\tau , v_{q_e}, \mu )
\right\rangle ,
\end{align*}
where $\zeta (\tau ,v_{q_{e}},\mu)
:=\mathbb{I}^{-1}((\operatorname{Exp}_{q_{e}}(\tau
v_{q_{e}}))\beta (\tau ,\mu )$. Since $\zeta (\tau ,v_{q_{e}},\mu
)$ is smooth in all variables also at $\tau =0$ by Remark
\ref{remark on smoothness}, it follows that
$\langle G_{1}(\tau ,[v_{q_{e}}],\mu ),\varsigma \rangle$ is a
smooth function of all its variables.
This expression at $\tau= 0 $ equals
\begin{align*}
&\langle G_{1}(0,[v_{q_{e}}],\mu ),\varsigma \rangle
=\mathbf{d}V(q_{e})
(\varsigma _{Q}(q_{e}))
-\frac{1}{2}\left\langle \zeta (0,
v_{q_e}, \mu),\; T_{q_e} \mathbb{I}
\left( \varsigma_Q(q_e)   \right)
\zeta (0 , v_{q_e}, \mu )
\right\rangle  \\
&\qquad =\mathbf{d}V(q_{e})
(\varsigma _{Q}(q_{e}))
-\frac{1}{2}\left\langle(\mathbb{I}(q_{e})[\zeta
(0,v_{q_{e}},\mu ),\varsigma ], \zeta (0,v_{q_{e}},\mu)
\right\rangle
- \frac{1}{2} \left\langle  \mathbb {I}(q_e)
\zeta(0,v_{q_{e}},\mu), [ \zeta(0,v_{q_{e}},\mu), \varsigma]
\right\rangle\\
&\qquad =\mathbf{d}V(q_{e})
(\varsigma _{Q}(q_{e}))
- \left\langle  \mathbb {I}(q_e)
\zeta(0,v_{q_{e}},\mu), [ \zeta(0,v_{q_{e}},\mu), \varsigma]
\right\rangle
\end{align*}
by \eqref{useful identity for I}. Since $V$ is $G$-invariant it
follows that  $\mathbf{d}V(q_{e})(\varsigma _{Q}(q_{e}))=0$.
Since $\zeta (0,v_{q_{e}},\mu ) = \xi(0,v_{q_{e}},\mu ) +
\eta_\mu \in \mathfrak{g}_{q_e} \oplus \mathfrak{k}_1 =
\mathfrak{t}$ (see Remark \ref{remark on smoothness}) it follows
that
$[\zeta (0,v_{q_{e}},\mu ),\varsigma ]\in [ \mathfrak{t},
\mathfrak{g}]$. By Proposition \ref{montaldi} (ii), we have
$\mathbb{I}(q_{e})\mathfrak{t} \subset [ \mathfrak{g},
\mathfrak{t}]^{\circ }$ and hence the second term above also
vanishes. Thus we get
$\langle G_{1}(0,[v_{q_{e}}],\mu ),\varsigma \rangle =0$ for any
$\varsigma \in \mathfrak{k}_2 $, that is, $G_{1}(0,[v_{q_{e}}],\mu
) = 0 $ which proves the theorem.
\end{proof}
\subsection{Bifurcating branches of relative equilibria}
Let $(Q,\langle\!\langle \cdot ,\cdot \rangle\!\rangle _{Q},V,G)$
be a simple mechanical $G$-system, with $G$ a compact Lie group
with the Lie algebra $\mathfrak{g}$. Let $ q_{e}\in Q$ be a
symmetric point whose isotropy group $G_{q_{e}}$ is contained in a
maximal torus $\mathbb{T}$ of $G$. Denote by $\mathfrak{t} \subset
\mathfrak{g}$ the Lie algebra of
$\mathbb{T}$.
Let $B\subset (\mathfrak{g}\cdot q_{e})^{\perp }$ be a
$G_{q_{e}}$--invariant open neighborhood of $0_{q_{e}}\in
(\mathfrak{g}\cdot q_{e})^{\perp }$ such that the exponential map
is injective on $B $ and for any $ q \in G\cdot
\operatorname{Exp}_{q_{e}}(B)$ the isotropy subgroup $G_{q} $ is
conjugate to a (not necessarily proper) subgroup of $G_{q_{e}}$.
Define the closed $G_{q_e}$--invariant subset
$\mathcal{Z}_{\mu^{0}}=:\{v_{q_{e}}\in B\cap
(T_{q_{e}}Q)_{\{e\}}\mid \det A=0\}$, where ${\mu^{0}}\in
\mathfrak{m}_{1}\oplus \mathfrak{m}_{2}$ is arbitrarily chosen and
the entries of the matrix $A $ are given in \eqref{A}.  Let $U
\subset [B\cap (T_{q_{e}}Q)_{\{e\}} \setminus
\mathcal{Z}_{\mu^{0}}]/G_{q_{e}} $ be open and consider the
functions $F$ and $G $ given in Theorems \ref{theorem about F} and
\ref{theorem about G}.
Define $G^{i}:I\times U\times (\mathfrak{m}_{1}\oplus
\mathfrak{m}_{2})\rightarrow \mathbb{R}$ by
\[
G^{i}(\tau,[v_{q_{e}}],\mu_{1}+\mu _{2}):=\langle
G(\tau,[v_{q_{e}}],\mu_{1}+\mu _{2}),\varsigma_{i}\rangle,
\]
where $\{\varsigma_{i} \mid  i=1,... ,\operatorname{dim}
\mathfrak{k_{2}}\}$ is a basis for ${\mathfrak k}_{2}$.
Choose $([v_{q_{e}}],\mu _{1}+\mu _{2})\in U\times
(\mathfrak{m}_{1}\oplus \mathfrak{m}_{2})$ such that
\[
\frac{\partial F}{\partial u}(0,[v_{q_{e}}],\mu _{1}+\mu _{2})=0,
\]
where the partial derivative is taken relative to the variable $u
\in U $. Define the matrix
\begin{equation*}
\Delta_{([v_{q_{e}}],\mu _{1},\mu _{2})} := \left[
\begin{array}{cc}
\frac{\partial ^{2}F}{\partial u^{2}}(0,[v_{q_{e}}],\mu _{1}+\mu
_{2}) & \frac{\partial^{2}F}{\partial \mu_{2}\partial
u}(0,[v_{q_{e}}],\mu _{1}+\mu _{2})\\
\frac{\partial G^{i}}{\partial u}(0,[v_{q_{e}}],\mu _{1}+\mu _{2})
& \frac{\partial G^{i}}{\partial \mu_{2}}(0,[v_{q_{e}}],\mu
_{1}+\mu _{2})
\end{array}
\right],
\end{equation*}
where the partial derivatives are evaluated at $\tau
=0, [v_{q_{e}}],\mu = \mu_1 + \mu_2$. Here $ \frac{\partial}
{\partial \mu _{2}}$ denotes the partial derivative with respect to
the $\mathfrak{m}_{2}$-component $\mu_2 $ of $\mu $.
In the framework and the notations introduced above we will state
and prove the main result of this paper. Let $\pi :TQ\rightarrow
(TQ)/G$ be the canonical projection and $\mathcal{R} _{e}:=\pi
(\mathfrak{t}\cdot q_{e})$.
\begin{theorem}
\label{principala}
Assume the following:
\begin{equation*}
{\bf (H)}\; \text{every } v_{q_{e}}\in \mathfrak{t}\cdot q_{e}
\text{ is a relative equilibrium.}
\end{equation*}
If there is a point $([v_{q_{e}}^{0}],\mu
_{1}^{0}+\mu _{2}^{0})\in U\times (\mathfrak{m}_{1}\oplus
\mathfrak{m}_{2})$ such that
\begin{eqnarray*}
&&
\begin{array}{cc}
1) & \frac{\partial F}{\partial u}(0,[v_{q_{e}}^{0}],\mu
_{1}^{0}+\mu _{2}^{0})=0,
\end{array}
\\
&&
\begin{array}{cc}
2) & G^{i}(0,[v_{q_{e}}^{0}],\mu _{1}^{0}+\mu _{2}^{0})=0
\end{array}
\\
&&
\begin{array}{cccc}
3) & \Delta _{([v_{q_{e}}^{0}],\mu _{1}^{0},\mu _{2}^{0})} & is &
non degenerate,
\end{array}
\end{eqnarray*}
then there exists a family of continuous curves
$\gamma_{_{({[v_{q_{e}}^{0}]},\mu _{1}^{0},\mu
_{2}^{0})}}^{\mu_{1}}:[0,1]\rightarrow (TQ)/G$
parameterized by
$\mu_{1}$ in a small neighborhood $\mathcal{V}_{0}$ of $\mu_{1}^{0}$
consisting of classes of relative equilibria with trivial
isotropy on $\gamma _{_{({[v_{q_{e}}^{0}]},\mu _{1}^{0},\mu
_{2}^{0})}}^{\mu_{1}}(0,1)$ satisfying
\begin{equation*}
\operatorname{Im}\gamma _{_{({[v_{q_{e}}^{0}]},\mu _{1}^{0}, \mu
_{2}^{0})}}^{\mu_{1}}\bigcap \mathcal{R} _{e}=\left\{\gamma
_{_{({[v_{q_{e}}^{0}]}, \mu _{1}^{0}, \mu
_{2}^{0})}}^{\mu_{1}}(0)\right\}
\end{equation*}
and $\gamma _{_{({[v_{q_{e}}^{0}]}, \mu _{1}^{0}, \mu
_{2}^{0})}}^{\mu_{1}}(0) = [ \zeta_Q(q_e)]$, where $\zeta =
\widehat{\mathbb {I}}(q_e)^{-1} \mu_1 \in \mathfrak{t}$.
For $\mu_{1},\mu_{1}^{\prime}\in \mathcal{V}_{0}$ with
$\mu_{1}\neq\mu_{1}^{\prime}$, where $\mathcal{V}_{0}$ is as
above, the above branches do not intersect, that is,
\begin{equation*}
\left\{ \gamma _{_{({[v_{q_{e}}^{0}]}, \mu _{1}^{0}, \mu
_{2}^{0})}}^{\mu_{1}}(\tau) \,\Big|\, {\tau\in [0,1]} \right\}
\bigcap
\left\{\gamma _{_{({[v_{q_{e}}^{0}]}, \mu _{1}^{0}, \mu
_{2}^{0})}}^{\mu_{1}^{\prime}}(\tau) \,\Big|\, {\tau\in
[0,1]}\right\} = \varnothing.
\end{equation*}
Suppose that $ ([v_{q_{e}}^{0}],\mu _{1}^{0},\mu
_{2}^{0}) \neq ([v_{q_{e}}^{1}],\mu _{1}^{1},\mu _{2}^{1})$.
(i) If $\mu_{1}^{0}\neq \mu_{1}^{1}$ then the families of relative
equilibria do not intersect, that is,
\begin{equation*}
\left\{\gamma_{_{({[v_{q_{e}}^{0}]}, \mu _{1}^{0},
\mu_{2}^{0})}}  ^{\mu_{1}}(\tau) \,\Big|\, {(\tau,\mu_{1})\in
[0,1]\times \mathcal{V}_{0}} \right\}
\bigcap
\left\{\gamma _{_{({[v_{q_{e}}^{1}]}, \mu_{1}^{1},
\mu_{2}^{1})}}^{\mu_{1}^{\prime}}(\tau) \,\Big|\,
{(\tau,\mu_{1}^{\prime})\in [0,1]\times \mathcal{V}_{1}}\right\}
=\varnothing,
\end{equation*}
where $\mathcal{V}_{0}$ and $\mathcal{V}_{1}$ are
two small neighborhoods of $\mu_{1}^{0}$ and $\mu_{1}^{1}$
respectively such that $\mathcal{V}_{0}\cap
\mathcal{V}_{1}=\varnothing$.
(ii) If $\mu_{1}^{0}=\mu_{1}^{1}=\overline \mu$ and
$[v_{q_{e}}^{0}]\neq[v_{q_{e}}^{1}]$ then $\gamma
_{_{({[v_{q_{e}}^{0}]},\overline\mu , \mu
_{2}^{0})}}^{\overline\mu}(0)=\gamma
_{_{({[v_{q_{e}}^{1}]},\overline\mu , \mu
_{2}^{1})}}^{\overline\mu}(0)$ and for $\tau>0$ we have
\begin{equation*}
\left\{\gamma
_{_{({[v_{q_{e}}^{0}]},\overline\mu , \mu
_{2}^{0})}}^{\overline\mu}(\tau) \,\Big|\, {\tau\in
(0,1]}\right\}
\bigcap \left\{\gamma
_{_{({[v_{q_{e}}^{1}]},\overline\mu , \mu
_{2}^{1})}}^{\overline\mu}(\tau) \,\Big|\, {\tau\in (0,1]}\right\}
=\varnothing.
\end{equation*}
\end{theorem}
\begin{proof}
Let $({[v_{q_{e}}^{0}]},\mu _{1}^{0}+\mu _{2}^{0})\in U\times
(\mathfrak{m}_{1}\oplus \mathfrak{m} _{2}) $ be such that the
conditions 1-3 hold. Because $\Delta _{({[v_{q_{e}}^{0}]},\mu
_{1}^{0}+\mu _{2}^{0})}$ is nondegenerate, we can apply the
implicit function theorem for the system $(\frac {\partial
F}{\partial u},G^{i})(\tau,[v_{q_{e}}],\mu _{1}+\mu _{2})=0$
around the point $(0,{[v_{q_{e}}^{0}]},\mu _{1}^{0}+\mu
_{2}^{0})$ and so we can find an open neighborhood $J\times
\mathcal{V}_{0}$ of the point $(0,\mu_{1}^{0})$ in $I\times
\mathfrak {m}_{1}$ and two functions $u:J\times
\mathcal{V}_{0}\rightarrow U$ and $\mu _{2}:J\times
\mathcal{V}_{0}\rightarrow \mathfrak{m}_{2}$ such that
$u(0,\mu _{1}^{0})={[v_{q_{e}}^{0}]}$, $\mu _{2}(0,\mu
_{1}^{0})=\mu _{2}^{0}$ and
\begin{eqnarray*}
&&
\begin{array}{cc}
i) & \frac{\partial F}{\partial u}(\tau ,u(\tau,\mu _{1} ),\mu
_{1}+\mu _{2}(\tau,\mu _{1}))=0
\end{array}
\\
&&
\begin{array}{cc}
ii) & G^{i}(\tau ,u(\tau,\mu _{1} ),\mu _{1}+\mu _{2}(\tau,\mu
_{1} ))=0.
\end{array}
\end{eqnarray*}
Therefore, from Theorems \ref{theorem about F} and \ref{theorem
about G} it follows that the relative equilibrium conditions of
Corollary \ref{criteriu} are both satisfied. Thus
we obtain the following family of branches of relative equilibria
$[(\overline{\sigma}(\tau\cdot
u(\tau,\mu_{1})),\beta(\tau,\mu_{1}+\mu_{2}(\tau,\mu_{1})))]_{G}$
parameterized by $\mu_{1} \in \mathcal{V}_{0}$. For $\tau>0$ the
isotropy subgroup is trivial and for
$\tau=0$ the corresponding points on the branches are
$[(\overline{\sigma}([0_{q_{e}}]),\mu_{1}]_{G}=[q_{e},\mu_{1}]_{G}$
which have the isotropy subgroup equal to $G_{q_{e}}$. This shows
that there are points in $\mathcal{R}_e $ from which there are
emerging branches of relative equilibria with broken trivial
symmetry.
Using now the correspondence given by Proposition \ref{map f} and
a rescaling of $\tau$ we obtain the desired family of continuous
curves $\gamma _{_{({[v_{q_{e}}^{0}]},\mu _{1}^{0},\mu
_{2}^{0})}}^{\mu_{1}}:[0,1]\rightarrow (TQ)/G$ parameterized by
$\mu_{1}$ in a small neighborhood $\mathcal{V}_{0}$ of $\mu_{1}^{0}$
consisting of classes of relative equilibria with trivial
isotropy on $\gamma _{_{({[v_{q_{e}}^{0}]},\mu _{1}^{0},\mu
_{2}^{0})}}^{\mu_{1}}(0,1)$ and such that
\begin{equation*}
\operatorname{Im}\gamma _{_{({[v_{q_{e}}^{0}]},\mu _{1}^{0},\mu
_{2}^{0})}}^{\mu_{1}}\bigcap \mathcal{R} _{e}=\{\gamma_{
_{({[v_{q_{e}}^{0}]},\mu _{1}^{0},\mu _{2}^{0})}}^ {\mu_1}(0)\}
\end{equation*}
and $\gamma_{_{({[v_{q_{e}}^{0}]},\mu _{1}^{0},\mu _{2}^{0})}}
^{\mu_1}(0) = [ \zeta_Q(q_e)]$, where $\zeta =
\widehat{\mathbb {I}}(q_e)^{-1} \mu_1$.
Equivalently, using the identification given by
(\ref{corespondenta}) and by Proposition \ref{map f} we obtain that
the branches of relative equilibria $\gamma
_{_{({[v_{q_{e}}^{0}]},\mu _{1}^{0},\mu
_{2}^{0})}}^{\mu_{1}}(\tau)\in (TQ)/G$ are identified with
$[\sigma(\tau\cdot
u(\tau,\mu_{1})),\beta(\tau,\mu_{1}+\mu_{2}(\tau,\mu_{1}))]
_{G_{q_{e}}}$. It is easy to see that for
$\mu_{1}\neq\mu_{1}^{\prime}$ we have that
$\beta(\tau,\mu_{1}+\mu_{2}(\tau,\mu_{1}))\neq\beta(\tau^{\prime},
\mu_{1}^{\prime}+\mu_{2}(\tau,\mu_{1}^{\prime}))$ for every
$\tau,\tau^{\prime}\in [0,1]$. Using now the fact that $G_{q_{e}}$
acts trivially on $\mathfrak{m}_{1}$ we obtain
\begin{equation*}
\left\{ \gamma _{_{({[v_{q_{e}}^{0}]},\mu _{1}^{0},\mu
_{2}^{0})}}^{\mu_{1}}(\tau) \,\Big| \, {\tau\in [0,1]} \right\}
\bigcap
\left\{\gamma _{_{({[v_{q_{e}}^{0}]},\mu _{1}^{0},\mu
_{2}^{0})}}^{\mu_{1}^{\prime}}(\tau) \,\Big| \, {\tau\in
[0,1]}\right\} = \varnothing.
\end{equation*}
In an analogous way, using the same argument we can prove $(i)$.
For $(ii)$ we start with two branches of relative equilibria,
$b_{1}(\tau,\overline\mu):=[\sigma(\tau\cdot
u(\tau,\overline\mu)),\beta(\tau,\overline\mu+\mu_{2}
(\tau,\overline\mu))]_{G_{q_{e}}}$ and
$b_{2}(\tau^{\prime},\overline\mu):=[\sigma(\tau^{\prime}\cdot
u^{\prime}(\tau^{\prime},\overline\mu)),\beta(\tau^{\prime},
\overline\mu+\mu_{2}(\tau,\overline\mu))]_{G_{q_{e}}}$. For
$\tau=\tau^{\prime}=0$ we have
$b_{1}(0,\overline\mu)=[0,\overline\mu]_{G_{q_{e}}}
=b_{2}(0,\overline\mu)$. We also have
$u(0,\overline\mu)=[v_{q_{e}}^{0}]\neq
[v_{q_{e}}^{1}]=u^{\prime}(0,\overline\mu)$ and so, from
the implicit  function theorem, we obtain
$u(\tau,\overline\mu)\neq u^{\prime}(\tau^{\prime},
\overline\mu)$ for $\tau,\tau^{\prime}>0$ small
enough.
Suppose that there exist  $\tau,\tau^{\prime}>0$ such that
$b_{1}(\tau,\overline\mu)=b_{2}(\tau^{\prime},\overline\mu)$. Then
using the triviality of the ${G_{q_{e}}}$-action on
$\mathfrak{m}_{0}$ we obtain that
$\tau^{2}\nu_{0}={\tau^{\prime}}^{2}\nu_{0}$ and consequently
$\tau=\tau^{\prime}$.  The conclusion of $(ii)$ follows now by
rescaling.
\end{proof}
\begin{remark}
\normalfont
We can have two particular forms for the rescaling $\beta $
according to special choices of the groups $G$ and  $G_{q_{e}}$,
respectively.
(a) If $G$ is a torus, then from the splitting $\mathfrak{g}=
\mathfrak{k}_{0}\oplus \mathfrak{k}_{1}\oplus \mathfrak{k}_{2}$, where $
\mathfrak{k}_{0}=\mathfrak{g}_{q_{e}}$, $\mathfrak{k}_{0}\oplus \mathfrak{k}
_{1}=\mathfrak{t}$, and
$\mathfrak{k}_{2}=[\mathfrak{g},\mathfrak{t}]$, we conclude
that $\mathfrak{k}_{2}=\{0\}$ (since $\mathfrak{g}= \mathfrak{t}$)
and consequently
$\mathfrak{m} _{2}=\{0\}$. In this case we will obtain the special
form for the rescaling $
\beta :I\times \mathfrak{m}_{1}\rightarrow \mathfrak{g}^{\ast }$,
$\beta (\tau ,\mu )=\mu +\tau ^{2}\nu _{0}$.
(b) If is $G_{q_{e}}$ a maximal torus in $G$, so $\mathfrak{g}_{q_e}
= \mathfrak{t}$, then the same splitting implies that
$\mathfrak{k}_{1}=\{0\}$ and consequently $\mathfrak{m}_{1}=\{0\}$.
In this case we will obtain the special form for the rescaling
$\beta :I\times \mathfrak{m} _{2}\rightarrow
\mathfrak{g}^{\ast }$, $\beta (\tau ,\mu )=\tau \mu +\tau ^{2}\nu
_{0}$.
\end{remark}

\section{Stability of the bifurcating branches of relative
equilibria}
\label{stability section}
In this section we shall study the stability of the
branches of relative equilibria found in the previous section.
We will do this by applying a result of Patrick \cite{patrick 
thesis} on $G_{\mu}$-stability to our situation. First we shortly
review this result.
\begin{definition}
Let $z_{e}$ be a relative equilibrium with velocity $\xi_{e}$ and
$J(z_{e})=\mu_{e}$. We say that $z_{e}$ is {\bfi formally
stable\/} if $\mathbf{d}^{2}(H-J^{\xi_{e}})(z_{e})|_
{T_{z_{e}}J^{-1}(\mu_{e})}$ is a positive or negative definite
quadratic form on some (and hence any) complement to ${\mathfrak
g}_{\mu_{e}}\cdot z_{e}$ in $T_{z_{e}}J^{-1}(\mu_{e})$.
\end{definition}
We have the following criteria for formal stability.
\begin{theorem}[Patrick, 1995]
\label{patrick}
Let $z_{e}\in T^{\ast }Q$ be a relative equilibrium with momentum
value $\mu _{e} \in \mathfrak{g}^\ast$ and base point $q_e \in
Q$. Assume that $\mathfrak{g}_{q_{e}}=\{0\}$. Then
$z_{e}$ is formally stable if and only if
$\mathbf{d}^{2}V_{\mu }(q_{e})$ is positive definite on one (and
hence any) complement $\mathfrak{g}_{\mu }\cdot q_{e}$ in
$T_{q_{e}}Q$.
\end{theorem}
To apply this theorem to our case in order to obtain the formal
stability of the relative equilibria on a bifurcating branch
we proceed as follows. First notice that if we fix $\mu
\in
\mathfrak{m}_{1}\oplus \mathfrak{ m}_{2}$ and $[v_{q_{e}}]\in U$
as in Theorem \ref{principala}, we obtain locally a branch
of relative equilibria with trivial isotropy bifurcating from
our initial set. More precisely, this branch starts at the
point 
\[
\left(\widehat{ \mathbb{I}}(q_{e})^{-1}\Pi_{1}\mu \right)_Q(q_e).
\]
The momentum values along this branch are  
$\beta (\tau ,\mu )$, and for $\tau
\neq 0$ the velocities have the expression $\mathbb{I}
(\operatorname{Exp} _{q_{e}}(\sigma(\tau
u(\tau,\mu_{1}))^{-1}\beta (\tau ,\mu )$. The base points of this
branch are
$\operatorname{Exp}_{q_{e}}(\sigma(\tau u(\tau,\mu_{1}))$.
Recall from Corollary \ref{criteriu} that we
introduced the notation $\overline{\sigma}  : =
\operatorname{Exp}{q_e} \circ \sigma $ that will be used below.
By the definition of $\beta (\tau ,\mu )$ we have 
$\mathfrak{g}_{\beta (\tau ,\mu )}=\mathfrak{t}$ for all $\tau $,
even for $\tau =0$. The base points for the
entire branch have no symmetry for $\tau>0 $ so we can
characterize the formal stability (in our case the $\mathbb{T}$
-stability) of the whole branch (locally) in terms of  Theorem
\ref{patrick}. We begin by giving sufficient conditions that
guarantee the $\mathbb{T}$-stability of the branch, since
$G_{\beta(\tau,
\mu) } = \mathbb{T} $. To do this, one needs to find conditions
that insure that for
$\tau
\neq 0 $ (where the amended  potential exists)
\[
\mathbf{d}^2 V_{\beta(\tau, \mu)}(
\overline{\sigma}(\tau u(\tau,\mu_{1}))|_{T_{[\tau
u(\tau,\mu_{1})]}
\overline{\sigma}(T_{[\tau u(\tau,\mu_{1})]}U) \oplus
(T_{\sigma ([\tau u(\tau,\mu_{1})])}\operatorname{Exp}_{q_{e}})
(\mathfrak{k} _{2}\cdot q_{e})}
\]
is positive definite. We do not know how to control the cross
terms of this quadratic form. This is why we shall work only with
Abelian groups $G $ since in that case the subspace
$\mathfrak{k}_2 = \{0\}$ and the second summand thus vanishes.
So, let $G $ be a torus $\mathbb{T}$. By Proposition
\ref{decomposition of the tangent space adapted to sigma bar} and
Theorem
\ref{theorem about F}, the second variation 
\[
\mathbf{d}^2 V_{\beta(\tau, \mu)}(
\overline{\sigma}(\tau u(\tau,\mu_{1}))|_{T_{[\tau
u(\tau,\mu_{1})]}
\overline{\sigma}(T_{[\tau u(\tau,\mu_{1})]}U)}
\]
coincides for $\tau\neq 0 $,  with the second variation 
\begin{equation}
\label{hessian of f one}
\mathbf{d}_{U}^2 F_1(\tau,
u(\tau,\mu_{1}), \mu_1 + \mu_2( \tau, \mu_1))|_{T_{[\tau
u(\tau,\mu_{1})]}U}
\end{equation}
of the auxiliary function $F_1 $, where $\mathbf{d}^2_{U}$ denotes the second variation relative to the
second variable in $F_1$.  But, unlike
$V_{\beta(\tau,\mu)}$, the function
$F_1$ is is defined even at $\tau = 0$. Recall from Theorem
\ref{theorem about F} that on  the bifurcating branch the amended
potential has the expression
\begin{equation*}
F_{1}(\tau ,u(\tau,\mu_{1}),\mu_1 + \mu_2( \tau, \mu_1)
)=F_{0}(\mu_1 + \mu_2( \tau, \mu_1))+\tau ^{2}F(\tau,
u(\tau,\mu_{1}),\mu_1 + 
\mu_2(\tau, \mu_1) ),
\end{equation*}
where $F_0$ is smooth on $\mathbf{J}_L (\mathfrak{g}\cdot
q_e) = \mathbb {I}(q_e)\mathfrak{g}$ and $F, F_1 $ are both 
smooth functions on $I\times U\times \mathbf{J}
_{L}(\mathfrak{g}\cdot q_{e})$, even around $\tau =0$. So, if the
second variation of $F $ at $(0, [v_{q_e}^0], \mu_1 ^0 + \mu_2
^0)$ is positive definite, then the quadratic form \eqref{hessian
of f one} will remain positive definite along the branch for $\tau
> 0 $ small. So we get the following result.
\begin{theorem}
Let $\mu_1 ^0 + \mu_2 ^0 \in \mathfrak{m}_{1}\oplus
\mathfrak{m}_{2}$ and
$[v_{q_{e}}^0]\in U$  be as in the Theorem \ref{principala} and
assume that
$\mathbf{d}^2_U F(0,[v_{q_{e}}^0],\mu_1^0 + \mu_2 ^0 )$ is
positive definite. Then the branch of relative
equilibria with no symmetry  which bifurcate form
$\left(\widehat{ \mathbb{I}}(q_{e})^{-1}\mu_1^0
\right)_Q(q_e)
$ will be $\mathbb{T}$-stable for $\tau>0$ small.
\end{theorem}
\medskip
A direct application of this criterion to the double spherical
pendulum recovers the stability result on the bifurcating branches
proved directly in \cite{ms}.
\medskip
\addcontentsline{toc}{section}{Acknowledgments}
\noindent\textbf{Acknowledgments.} We would like to thank J.
Montaldi for telling us the result of Proposition \ref{montaldi}
and for many discussions that have influenced our presentation
and clarified various points during the writing of this paper.
Conversations with A. Hern\'andez and J. Marsden are also
gratefully acknowledged. The third and fourth authors were
partially supported by the European Commission and the Swiss
Federal Government through funding for the Research Training Network
\emph{Mechanics and Symmetry in Europe} (MASIE). The first and
third author thank the Swiss National Science Foundation for
partial support.

\medskip
\noindent {\sc P. Birtea and M. Puta}, \\
Departamentul de Matematic\u a, Universitatea de Vest,
RO--1900 Timi\c soara, Romania.\\
Email: {\sf birtea@geometry.uvt.ro, puta@geometry.uvt.ro}\\
\noindent {\sc T.S. Ratiu and R\u azvan Tudoran},\\
Section de math\'ematiques,
{\'E}cole Polytechnique F{\'e}d{\'e}rale de Lausanne. CH--1015 Lausanne.
Switzerland.\\
Email: {\sf tudor.ratiu@epfl.ch, razvan.tudoran@epfl.ch}
\end{document}